\documentclass[review,11pt]{elsarticle}

\usepackage{fullpage}
\usepackage{amsmath,amssymb}
\usepackage[toc]{appendix}
\usepackage[T1]{fontenc}
\usepackage{textcomp}
\usepackage{graphicx,epstopdf}
\usepackage[position=top,labelfont=normalfont,textfont=normalfont,singlelinecheck=off,justification=raggedright]{subfig}
\usepackage{floatrow}
\usepackage[dvipsnames]{xcolor}
\usepackage{setspace}
\usepackage{fancyhdr}
\usepackage{enumerate,etaremune}
\usepackage{tabularx,multirow}
\usepackage{framed}
\usepackage{hhline}
\usepackage[pdftex,unicode,bookmarks=false]{hyperref}
\usepackage{url}
\hypersetup{
  colorlinks,
  citecolor=Blue, linkcolor=Blue, urlcolor=Blue
}
\usepackage{tikz}
\usetikzlibrary{shapes.geometric}
\usetikzlibrary{shapes,arrows}
\usetikzlibrary{plotmarks}

\usetikzlibrary{calc}
\usepackage{algorithm}
\usepackage[noend]{algpseudocode}

\newcommand{\eg}{\textit{e.g.}}
\newcommand{\ie}{\textit{i.e.}}
\newcommand{\cf}{\textit{cf.}}
\newcommand{\etal}{\textit{et al.}}
\newcommand{\tensor}[1]{\ensuremath{\boldsymbol{#1}}}
\newcommand{\transpose}[1]{{#1}^{\mathsf T}}

\newcommand{\od}{\mathrm{d}}
\newcommand{\pd}{\partial}
\newcommand{\Del}{\mathrm{\Delta}}

\newcommand{\inv}{\mathrm{inv}}
\newcommand{\vis}{\mathrm{vis}}
\DeclareMathOperator{\grad}{\nabla}
\DeclareMathOperator{\diver}{\grad\cdot}
\DeclareMathOperator{\symgrad}{\nabla^{s}}

\DeclareMathOperator{\tr}{tr}

\newsavebox{\dotbox}

\newcommand{\revised}[1]{{\color{black} #1}}

\newcolumntype{L}[1]{>{\raggedright\let\newline\\arraybackslash\hspace{0pt}}m{#1}}
\newcolumntype{C}[1]{>{\centering\let\newline\\arraybackslash\hspace{0pt}}m{#1}}
\newcolumntype{R}[1]{>{\raggedleft\let\newline\\arraybackslash\hspace{0pt}}m{#1}}

\allowdisplaybreaks

\biboptions{sort&compress,square,comma,numbers}
\AtBeginDocument{\hypersetup{citecolor=MidnightBlue,linkcolor=MidnightBlue,urlcolor=MidnightBlue}}

\usepackage{etoolbox}
\makeatletter
\patchcmd{\ps@pprintTitle}% <cmd>
{Preprint submitted to}% <search>
{~}% <replace>
{}{}% <succes><failure>
\makeatother
\begin{document}

\begin{frontmatter}

\title{Stabilized material point methods for coupled large deformation and fluid flow in porous materials}

\author[HKU]{Yidong Zhao}
\author[HKU]{Jinhyun Choo\corref{corr}}
\ead{jchoo@hku.hk}

\cortext[corr]{Corresponding Author}

\address[HKU]{Department of Civil Engineering, The University of Hong Kong, Pokfulam, Hong Kong}

\journal{~}

\begin{abstract}
The material point method (MPM) has been increasingly used for the simulation of large-deformation processes in fluid-infiltrated porous materials.
For undrained poromechanical problems, however, standard MPMs are numerically unstable because they use low-order interpolation functions that violate the inf--sup stability condition.
In this work, we develop stabilized MPM formulations for dynamic and quasi-static poromechanics that permit the use of standard low-order interpolation functions notwithstanding the drainage condition.
For the stabilization of both dynamic and quasi-static formulations, we utilize the polynomial pressure projection method whereby a stabilization term is augmented to the balance of mass.
The stabilization term can be implemented with both the original and generalized interpolation material point (GIMP) methods, and it is compatible with existing time-integration methods.
Here we use fully-implicit methods for both dynamic and quasi-static poromechanical problems, aided by a block-preconditioned Newton--Krylov solver.
The stabilized MPMs are verified and investigated through several numerical examples under dynamic and quasi-static conditions.
Results show that the proposed MPM formulations allow standard low-order interpolation functions to be used for both the solid displacement and pore pressure fields of poromechanical formulations, from undrained to drained conditions, and from dynamic to quasi-static conditions.
\end{abstract}


\end{frontmatter}

% SECTION 1
% ------------------------------------------------------------------------------
\section{Introduction}
Large deformations in fluid-infiltrated porous materials are central to many important problems in science and engineering.
Notable examples are landslides, levee failures, and ground collapses, all of which continuously threaten and damage the built environment~\cite{Orlandini2015,Iverson2015,Borja2016b,Dey2016,Baer2018}.
Besides, soft porous materials in biophysics and materials science, such as tissues and hydrogels, often undergo very large deformations coupled with the flow of the pore fluid~\cite{Franceschini2006,Bertrand2016,Dortdivanlioglu2019}.
Therefore, the numerical modeling of large-deformation problems in porous materials has been a significant challenge for researchers in many disciplines alike.

A majority of numerical simulations of coupled poromechanics have employed mesh-based Lagrangian methods that trace solid material points (\eg~\cite{Borja1998,Borja1998a,Sanavia2002,Li2004,Uzuoka2012,Borja2016,Krischok2016,Na2017,Choo2018d,Dortdivanlioglu2018}).
These Lagrangian methods can straightforwardly incorporate loading and deformation history in the solid material, which is a crucial advantage over Eulerian methods because many porous materials such as soils show strong nonlinearity and history-dependence.
However, a purely Lagrangian method loses its capability when deformation gives rise to significant mesh distortion.
A successive remeshing technique may be used to overcome this drawback, but its application can be quite cumbersome especially when a number of history-dependent variables should be projected to the new mesh.
For these reasons, several numerical methods have been developed and advanced to simulate large deformations in a porous solid without mesh distortion while tracking the material's loading and deformation history efficiently (see~\cite{Chen2017} for a recent review of these methods).

The material point method (MPM) is a hybrid Lagrangian--Eulerian method initially proposed by Sulsky and coworkers~\cite{Sulsky1994,Sulsky1995} as a solid mechanics version of the fluid implicit particle (FLIP) method~\cite{Brackbill1986} which itself is an offspring of the particle-in-cell (PIC) method~\cite{Harlow1964}.
In essence, the MPM can be regarded as a variant of the finite element method (FEM) whereby quadrature points (material points) are allowed to move freely within a computational mesh.
The MPM's mesh is independent of the domain geometry and serves as a background mesh to calculate the weak form---like a mesh in an Eulerian method.
This feature makes the MPM immune to mesh distortion.
At the same time, because the MPM explicitly uses and traces material particles---like in a Lagrangian method---it is naturally able to incorporate the history-dependent solid behavior.
Thanks to these attractive features, the MPM has been gaining increasing popularity for the simulation of large-deformation processes in a variety of materials including fluid-infiltrated porous media~\cite{Zhang2009,Mackenzie-Helnwein2010,Abe2014,Bandara2015,Yerro2015,Bandara2016,Soga2016,Tampubolon2017,Liu2017,Pinyol2018}.

When applied to coupled poromechanics, MPM formulations face the challenge of how to address the inf--sup stability condition in the limit of undrained deformation/incompressible flow.
In the undrained limit, there is no relative flow between the pore fluid and the solid matrix, and so the pore pressure field acts as an incompressibility constraint of the solid deformation.
Then, because the MPM uses an interpolation scheme essentially the same as mixed finite elements, it is subjected to the same inf--sup stability condition for mixed finite elements for constrained problems~\cite{Brezzi1974,Brezzi1990,Bathe2001}.
In the context of poromechanics, the satisfaction of this condition requires the solid displacement field to be interpolated by higher-order functions than the pore pressure field (\ie~by the Taylor--Hood elements~\cite{Taylor1973}).

The stability requirement is particularly demanding for the MPM because all standard MPM formulations use low-order (linear) interpolation functions.
The MPM's dominant use of low-order interpolation functions is mainly due to computational efficiency.
In the MPM, the number of material points is often larger than the number of quadrature points in the standard FEM, to integrate weak forms with sufficient accuracy even when the material points are placed far from the optimal quadrature locations.
As such, if higher-order interpolation functions are employed, the number of material points should be dramatically increased for numerical integration.
For this reason, it will be impractical to use higher-order interpolation functions for the MPM, especially when iterations have to be performed in many material points due to nonlinear material behavior.

In the MPM literature, two types of approaches have been proposed to attain numerical stability without using higher-order interpolation functions.
The first type of approach, used in Jassim \etal~\cite{Jassim2013}, is to apply a fractional time stepping method originally developed for mixed finite elements for dynamic poromechanical problems~\cite{Pastor2000,Huang2001,Li2003}.
This method performs well for dynamic problems in the undrained limit, but it is conditionally stable in time and does not cover problems in drained conditions.
Moreover, as this method uses a particular type of time stepping algorithm, it requires significant efforts for modifying the existing MPM code.
The second type of approach, which was introduced in more recent works~\cite{Abe2014,Bandara2015}, is a reduced integration scheme that evaluates the gradients of pore pressure shape functions at cell centers instead of cell nodes.
Although this method has been demonstrated to improve numerical stability, it is still unable to completely resolve the stability problem in the undrained limit~\cite{Bandara2015}.
Note also that reduced integration in low-order elements should be carefully performed to avoid spurious energy modes.
Therefore, a more robust and versatile method is desired for addressing the stability problem in the MPM for coupled poromechanics.

In this paper, we develop stabilized MPM formulations for dynamic and quasi-static poromechanics that permit the use of standard low-order interpolation functions throughout the entire range of drainage condition.
This development draws on the polynomial pressure projection (PPP) method, which was initially developed for stabilization of mixed finite elements for the Stokes flow~\cite{Bochev2006a} and then extended to mixed finite elements for coupled poromechanics under quasi-static~\cite{White2008,Preisig2011,Sun2013a,Choo2015,Sun2015,Choo2019} and dynamic conditions~\cite{Monforte2019}.
This stabilization method has two key advantages over the fractional stepping method.
First, it can be applied to both quasi-static and dynamic problems, which is a highly desirable feature because fluid-infiltrated porous media can deform in a quasi-static manner (\eg~soil consolidation problems) as well as in a dynamic way (\eg~debris flow and earthquake problems).
Second, it can be easily implemented in the existing MPM code because it preserves the existing algorithms for time integration and spatial discretization.
To our knowledge, this is the first work that applies this type of stabilization method to the MPM for coupled poromechanics.

The paper is organized as follows.
Section~\ref{sec:coupled-formulation} briefly describes a continuum formulation for a coupled poromechanical problem at large strains.
We focus our description on a dynamic poromechanical formulation because it can be reduced to a quasi-static formulation when inertial effects are absent.
In Section~\ref{sec:discretization}, we describe MPM discretization of the poromechanical formulations under dynamic and quasi-static conditions.
In doing so, we present formulations and algorithms for a fully-implicit MPM for coupled poromechanics, for the first time to our knowledge.
In Section~\ref{sec:stabilization}, we develop stabilized versions of the dynamic and quasi-static poromechanical formulations in which stabilization terms are augmented to allow for the use of standard low-order interpolation functions even in the undrained limit.
In Section~\ref{sec:examples}, we verify and investigate the stabilized MPM formulations through several numerical examples under dynamic and quasi-static conditions.
Section~\ref{sec:closure} concludes the work.

% SECTION 2
% ------------------------------------------------------------------------------
\section{Coupled poromechanical formulation}
\label{sec:coupled-formulation}
This section presents a continuum formulation for a coupled poromechanical problem at large strains.
As both quasi-static and dynamic conditions are common in large-deformation poromechanics, we will mainly describe a more general dynamic formulation with an explanation about how it can boil down to a quasi-static formulation.
For brevity, we will write the essence of the formulation only and refer to the relevant literature for its thorough derivation~\cite{Bowen1982,Borja1995a,Li2004}.
Since the focus of this work is on numerical instability due to the incompressibility constraint in undrained limit,
we will assume that both the solid and fluid phases are strictly incompressible.
This incompressibility assumption will give rise to undrained conditions that are most challenging for stable solutions.
Note that the solid matrix (skeleton) can still be compressible depending on the constitutive law.
We will also assume that thermal, chemical, and other complicated effects are absent.
This type of formulation is commonly used for water-saturated soils in geotechnical engineering, among others.

\subsection{Mixture representation and kinematics}
Using mixture theory we idealize a fluid-saturated porous material as a two-phase mixture of solid and fluid.
Throughout this paper, we will use indices $s$ and $f$ to denote quantities pertaining to the solid and fluid phases, respectively.
These indices will be used as subscripts for intrinsic properties of the phases, and as superscripts for average or partial properties of the phases.
For this mixture, we define the volume fractions of the solid and fluid phases as $\phi^{s}$ and $\phi^{f}$, respectively, which satisfy $\phi^{s}+\phi^{f}=1$.
Similarly, we define the partial mass densities of the solid and fluid phases as $\rho^{s} := \phi^{s}\rho_{s}$ and $\rho^{f} := \phi^{f}\rho_{f}$, respectively, where $\rho_{s}$ and $\rho_{f}$ are the intrinsic mass densities of the solid and fluid phases.
Accordingly, the average mass density of the mixture is given by $\rho := \rho^{s} + \rho^{f}$.
Note that while $\rho_{s}$ and $\rho_{f}$ are constant due to the assumption of incompressible solid and fluid phases,
$\rho^{s}$ and $\rho^{f}$ (and thus $\rho$) are evolving by the volumetric deformation of the solid matrix.

For describing the kinematics of the mixture, we use a Lagrangian approach that follows the motion of the solid matrix.
The motion of the solid matrix by $\tensor{\varphi}(\tensor{X},t)$, where $\tensor{X}$ is the position vector of a solid material point $X$ in the reference configuration and $t$ is time.
Accordingly, the displacement vector of the solid material point is given by $\tensor{u}(\tensor{X},t) := \tensor{\varphi}(\tensor{X},t) - \tensor{X}$.
The deformation gradient of the solid motion is then defined as
\begin{align}
    \tensor{F} := \frac{\pd \tensor{\varphi}}{\pd \tensor{X}}\,.
\end{align}
The Jacobian is defined as
\begin{align}
    J := \det \tensor{F} = \od v / \od V > 0\,,
\end{align}
where $\od V$ and $\od v$ are the differential volumes of the porous material in the reference and current configurations, respectively.

The material time derivative following the solid motion is given by
\begin{align}
    \frac{\od (\circ)}{\od t} &:= \frac{\pd (\circ)}{\pd t} + \grad(\circ)\cdot\tensor{v}\,,
\end{align}
where $\grad$ is the gradient operator evaluated with respect to the current configuration, and $\tensor{v}$ is the intrinsic velocity of the solid.
Then the material time derivative following the fluid motion can be expressed as
\begin{align}
    \frac{\od^{f} (\circ)}{\od t} &:= \frac{\od (\circ)}{\od t} + \grad(\circ)\cdot\tilde{\tensor{v}}\,,\quad
    \tilde{\tensor{v}} := \tensor{v}_{f} - \tensor{v}\,,
\end{align}
where $\tensor{v}_{f}$ is the intrinsic velocity of the fluid, and $\tilde{\tensor{v}}$ is the relative velocity of the fluid.
Similarly, we also define the material acceleration vectors of the solid and fluid phases as
\begin{align}
    \tensor{a} := \frac{\od \tensor{v}}{\od t}\,, \quad
    \tensor{a}_{f} := \frac{\od^{f} \tensor{v}_{f}}{\od t}\,, \quad
\end{align}
respectively.
In the following, we will simply denote $\od(\circ)/\od t$ by the dot, \eg~$\tensor{v}=\dot{\tensor{u}}$ and $\tensor{a}=\dot{\tensor{v}}=\ddot{\tensor{u}}$.

Hereafter, we will assume that $\tensor{a}$ equals $\tensor{a}_{f}$, i.e. the fluid phase has no relative acceleration to the solid phase.
This assumption, which was introduced in many previous works (\eg~\cite{Zienkiewicz1984,Li2003,Li2004,Zhang2009,Uzuoka2012,Abe2014}), is reasonable for porous media under relatively low-frequency loading whereby mixing or separation of the solid and fluid phases are not significant.
Notably, Abe \etal~\cite{Abe2014} have shown that, even with this assumption, a poromechanical model can reproduce the dynamic failure pattern of a levee as observed from experiments.
From a numerical point of view, the upshot of this assumption is that it removes the need to discretize the fluid velocity vector field, which can save significant computational efforts for multi-dimensional problems.

\subsection{Governing equations}
The governing equations of the coupled poromechanical problem come from two balance laws: (a) the balance of linear momentum, and (b) the balance of mass.
Because the MPM uses an updated Lagrangian approach, we shall write these balance equations in the current configuration.
The balance of linear momentum of the mixture can be expressed as
\begin{align}
    \diver{\tensor{\sigma}} + \rho\tensor{g} = \rho\tensor{a}\,,
\end{align}
where $\diver$ is the divergence operator evaluated with respect to the current configuration,
$\tensor{\sigma}$ is the total Cauchy stress tensor in the mixture,
and $\tensor{g}$ is the gravitational \revised{acceleration} vector.
The balance of mass can be written as
\begin{align}
    \diver\tensor{v} + \diver{\tensor{q}} = 0\,,
\end{align}
where $\tensor{q}:=\phi^{f}\tilde{\tensor{v}}$ is the Darcy flux vector.

Constitutive models should be introduced to complete the formulation.
Without loss of generality, we shall assume isotropy for solid deformation and fluid flow properties.
For modeling the deformation of the solid matrix, we will adopt Terzaghi's effective stress theory, which can be written as
\begin{align}
    \tensor{\sigma} = \tensor{\sigma}' - p\tensor{1}\,,
\end{align}
where $\tensor{\sigma}'$ the effective Cauchy stress tensor,
$p$ is the pore pressure,
and $\tensor{1}$ is the rank-two identity tensor.
To consider viscous damping effects under dynamic loading, the effective stress tensor is considered the summation of an inviscid stress, $\tensor{\sigma}'_{\rm inv}$, and a viscous part, $\tensor{\sigma}'_{\rm vis}$, \ie
\begin{align}
    \tensor{\sigma}' = \tensor{\sigma}'_{\rm inv} + \tensor{\sigma}'_{\rm vis}\,.
\end{align}
The inviscid stress is modeled by rate-independent isotropic hyperelasticity, given by
\begin{align}
    \tensor{\sigma}'_{\inv} = \frac{1}{J}\left(2\tensor{b}\cdot\frac{\pd \varPsi(\tensor{b})}{\pd \tensor{b}}\right)\,,
\end{align}
where $\tensor{b} := \tensor{F}\cdot\transpose{\tensor{F}}$ is the left Cauchy--Green tensor
and $\varPsi$ is the strain energy density function.
Without loss of generality, in this work we use Hencky elasticity of which strain energy density function can be written as
\begin{align}
  \varPsi(\tensor{b}) := \frac{\lambda}{2}(\ln J)^{2} + G\tr\left(\frac{1}{2}\ln \tensor{b}\right)^{2}\,,
\end{align}
where $\lambda$ and $G$ are the Lam\'{e} parameters.
These parameters can be converted into an equivalent set of two elasticity parameters, such as the bulk modulus $K$ and Poisson's ratio $\nu$.
Next, adopting a Kelvin-type model, we calculate the viscous stress as~\cite{Simo1998,Uzuoka2012}
\begin{align}
    \tensor{\sigma}'_{\vis} = \alpha_{\vis}\,\tensor{c}:\symgrad\tensor{v}\,,
\end{align}
where $\alpha_{\vis}$ is the damping parameter,
$\tensor{c}$ is the rank-four spatial tangential tensor,
and $\symgrad$ denotes the symmetric gradient operator evaluated with respect to the current configuration.

As for fluid flow in porous media, we assume that Darcy's law holds and write
\begin{align}
    \tensor{q} = \revised{-}\kappa\,[\grad p - \rho_{f}(\tensor{g}-\tensor{a})]\,,
    \label{eq:darcy-law}
\end{align}
where $\kappa:=k/\mu_{f}$ is the mobility, where $k$ is the intrinsic permeability and $\mu_{f}$ is the dynamic viscosity of the fluid.
The permeability may be assumed to be constant under small to moderate deformation, but it has to be related to the material's porosity under sufficiently large deformation.
For the permeability--porosity relation, here we use the Kozeny--Carman equation, given by
\begin{align}
    k = k_{0}\left(\frac{(1 - \phi_{0})^{2}}{\phi_{0}^{3}}\right)\left(\frac{\phi^{3}}{(1 - \phi)^2}\right)\,,
    \label{eq:kozeny-carman}
\end{align}
where $\phi$ denotes the porosity, and the subscript $0$ is used to denote the initial values of the permeability and the porosity.

\subsection{Simplified form for quasi-static problems}
Inertial effects can be neglected in many large-deformation processes in fluid-saturated porous materials.
Examples include large-strain consolidation in soil~\cite{Gibson1967,Gibson1981}.
For such problems, the foregoing formulation can be simplified by setting
\begin{align}
    \tensor{a} = \tensor{a}_{f} = \revised{\tensor{0}} \quad \text{and} \quad
    \alpha_{\rm vis} = 0\,.
    \label{eq:quasi-static-condition}
\end{align}
Then, we recover the standard quasi-static formulation for coupled poromechanics at large strains~\cite{Borja1995a,Borja1998a}.
Because the quasi-static formulation is also important in poromechanics,
its stabilized MPM formulation will also be developed later in this work.
For brevity, in the sequel we will mainly present a dynamic formulation with descriptions of differences between dynamic and quasi-static formulations other than Eq.~\eqref{eq:quasi-static-condition}.

\subsection{Strong form}
Let us denote by $\Omega\in\mathbb{R}^{d}$ the $d$-dimensional domain in the current configuration and by $\partial\Omega$ the domain's boundary.
The boundary is decomposed into
displacement (Dirichlet) and traction (Neumann) boundaries, $\partial_{u}\Omega$ and $\partial_{t}\Omega$,
and pressure (Dirichlet) and flux (Neumann) boundaries, $\partial_{p}\Omega$ and $\partial_{q}\Omega$.
The boundary decomposition should satisfy
$\partial\Omega = \overline{\partial_{u}\Omega\cup\partial_{t}\Omega}
= \overline{\partial_{u}\Omega\cup\partial_{t}\Omega}$
and $\emptyset = \partial_{u}\Omega\cap\partial_{t}\Omega
= \partial_{p}\Omega\cap\partial_{q}\Omega$.
The time interval of the problem is denoted by $\mathcal{T}:=(0,T]$ with $T > 0$.

The strong form of an initial--boundary-value problem can be stated as follows:
Given $\hat{\tensor{u}}$, $\hat{\tensor{t}}$, $\hat{p}$, $\hat{q}$, $\tensor{u}_{0}$, and $p_{0}$,
find $\tensor{u}$ and $p$ such that
\begin{align}
    \diver\tensor{\sigma} + \rho\tensor{g} = \rho\tensor{a} \quad
    &\text{in} \quad \Omega \times \mathcal{T}\,,
    \label{eq:momentum-balance-strong-form}\\
    \diver\tensor{v} + \diver\tensor{q} = 0 \quad
    &\text{in} \quad \Omega \times \mathcal{T}\,,
    \label{eq:mass-balance-strong-form}
\end{align}
subjected to the boundary conditions ($\tensor{n}$ is the unit outward normal vector in the current configuration)
\begin{align}
    \tensor{u} = \hat{\tensor{u}} \quad
    &\text{on}\quad \partial_{u}\Omega \times \mathcal{T}\,, \\
    \tensor{n}\cdot\tensor{\sigma} = \hat{\tensor{t}} \quad
    &\text{on}\quad \partial_{t}\Omega\times \mathcal{T}\,,  \\
    p = \hat{p} \quad
    &\text{on}\quad \partial_{p}\Omega\times \mathcal{T}\,,  \\
    -\tensor{n}\cdot\tensor{q} = \hat{q} \quad
    &\text{on}\quad \partial_{q}\Omega\times \mathcal{T}\,,
\end{align}
and the initial conditions $\tensor{u}(\tensor{X},0) = \tensor{u}_{0}(\tensor{X})$ and $p(\tensor{X},0) = p_{0}(\tensor{X})$.

% SECTION 3
% ------------------------------------------------------------------------------
\section{Material point method discretization}
\label{sec:discretization}
This section describes MPM discretization of the coupled poromechanical formulation presented above.
For both dynamic and quasi-static problems, we use fully-implicit methods in conjunction with a block-preconditioned Newton--Krylov solver.
To our knowledge, the dynamic and quasi-static formulations presented herein are the first fully-implicit MPM for coupled poromechanics.

\subsection{Spatial discretization}
\begin{figure}[htbp]
  \centering
  \includegraphics[width=0.75\textwidth]{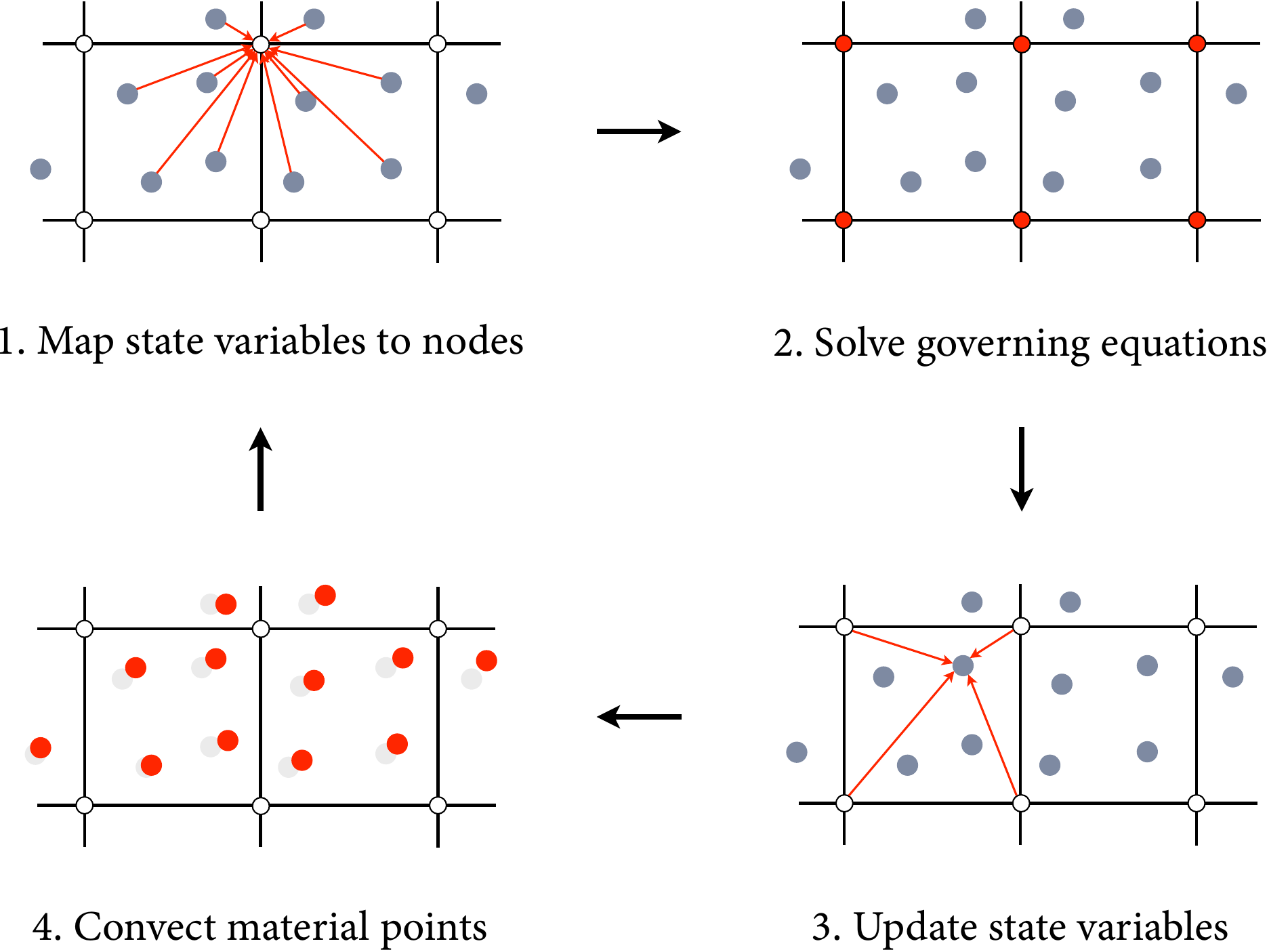}
  \caption{Illustration of the standard MPM procedure. A time (load) step consists of Stages 1--4. Once Stage 4 is finished, the next time step starts from Stage 1.}
  \label{fig:mpm-procedure}
\end{figure}

Figure~\ref{fig:mpm-procedure} illustrates the MPM procedure during a given time (load) increment.
The procedure may be divided into four stages, which can be briefly explained as follows:
% \begin{enumerate}[{Stage} 1]\itemsep=0pt
\begin{enumerate}[{Stage} 1.]\itemsep=0pt
  \item Map the state variables of material points (particles) to nodes in the background mesh.
  \item Solve discrete variational equations. The degrees of freedom are at the nodes, and the variational equations are integrated at the material points.
  \item Update the state variables of material points using solutions mapped from the nodes.
  \item Convect the material points according to the new variables. Reset the background grid if necessary.
\end{enumerate}
Here, Stage 2 is very similar to the standard FEM except that numerical integrations in the MPM are performed at material points, which move after each load step, instead of fixed quadrature points in the FEM.

In Stages 1 and 3, which are new to the MPM, state variables are mapped between material points and nodes.
For mapping an arbitrary variable $f$ from material points to nodes, the most widely used method in the MPM literature can be written as
\begin{align}
    f_{i} = \sum_{mp=1}^{n_{mp}}S_{i,mp}m_{mp}f_{mp}/M_{i}\,, \quad
    M_{i} := \sum_{mp=1}^{n_{mp}} S_{i,mp}m_{mp}\,.
    \label{eq:map-from-particles-to-nodes}
\end{align}
Here, subscripts $i$ and $mp$ are indices for nodes and material points, respectively, $n_{mp}$ is the number of material points in elements surrounding nodes, $m_{mp}$ is the mass of material point $mp$, and $S_{i,mp}$ is the MPM's basis function associated with node $i$ evaluated at the position of material point $mp$.
The basis functions are the same as the standard finite element shape functions in the original MPM, and they can be modified if desired.
Regardless, it is standard to use linear functions as the basis functions.
The above mapping method, which is standard in the MPM literature, uses a lumped mass matrix $M_{i}$ in Eq.~\eqref{eq:map-from-particles-to-nodes}.
However, Burgess \etal~\cite{Burgess1992} have shown that the use of the lumped mass matrix can lead to excessive numerical dissipation for a cyclic loading problem. They have proposed the following mapping method
\begin{align}
    \sum_{j=1}^{n_{d}} M_{ij}f_{j} = \sum_{mp=1}^{n_{mp}}S_{i,mp}m_{mp}f_{mp}\,, \quad
    M_{ij} := \sum_{mp=1}^{n_{mp}} S_{i,mp}S_{j,mp}m_{mp}\,,
    \label{eq:map-from-particles-to-nodes-consistent}
\end{align}
where $M_{ij}$ is the so-called consistent mass matrix, with $i$ and $j$ denoting node indices and $n_{d}$ the number of nodes surrounding the material point.
While this mapping improves accuracy for certain problems, it requires one to solve Eq.~\eqref{eq:map-from-particles-to-nodes-consistent} additionally.
In this work, we use the standard mapping method~\eqref{eq:map-from-particles-to-nodes} by default, adopting the consistent mass method~\eqref{eq:map-from-particles-to-nodes-consistent} for cyclic loading problems only.
The mapping from nodes to material points is conducted using the same basis functions
\begin{align}
    f_{mp} = \sum_{i=1}^{n_{d}} S_{i,mp}f_{i}\,. \quad
    \label{eq:map-from-nodes-to-particles}
\end{align}
More details of these MPM procedures are well described in the literature (see~\cite{Charlton2017,Sinaie2017,Zhang2017} for recent examples) so we shall focus on poromechanics-specific formulations in the following.

To begin the discretization procedure, we develop the variational form of the governing equations.
Let $\tensor{\eta}$ and $\psi$ denote the variations of $\tensor{u}$ and $p$, respectively.
Then the variational form of the linear momentum balance in the current configuration can be developed as
\begin{align}
    - \int_{\Omega} \symgrad\tensor{\eta} : \tensor{\sigma}'\, \od v
    - \int_{\Omega} (\diver\tensor{\eta}) p\, \od v
    + \int_{\Omega} \tensor{\eta} \cdot \rho (\tensor{g} - \tensor{a})\, \od v
    + \int_{\partial\Omega} \tensor{\eta} \cdot \hat{\tensor{t}}\,\od a = 0\,,
    \label{eq:momentum-weak-form-ver1}
\end{align}
and the variational form of the mass balance as
\begin{align}
    \int_{\Omega} \psi \diver{\tensor{v}}\, \od v
    - \int_{\Omega} \grad{\psi} \cdot \tensor{q}\, \od v
    + \int_{\partial\Omega} \psi \hat{q}\,\od a = 0\,.
    \label{eq:mass-weak-form-ver1}
\end{align}

Next, we use the Galerkin method to discretize the variational equations in space.
In the original MPM, the basis functions for spatial interpolation are identical to the standard FE shape functions denoted by $N$ above.
However, it is well known that this original MPM may suffer from oscillations when materials are crossing cell boundaries in the background grid.
A popular remedy for this cell-crossing instability is the generalized interpolation material point (GIMP) method~\cite{Bardenhagen2004}.
In a local 1D coordinate $\xi$, the GIMP method sets the basis function as
\begin{align}
  S_{i,mp} := \frac{1}{V_{mp}}\int_{\Omega_{mp}\cap\Omega} \chi_{mp}(\xi)N_{i}(\xi)\,\od\xi\,,
  \label{eq:gimp-basis}
\end{align}
where $V_{mp}$ and $\Omega_{mp}$ are the volume and the influence domain of material point $mp$, respectively, and $\chi_{mp}$ is the indicator function defined as
\begin{align}
  \chi_{mp}(\xi) =
  \left\{\begin{array}{ll}
  1\,, & \text{if}\;\, \xi\cap\Omega_{mp}\,, \\
  0\,, & \text{otherwise}\,.
  \end{array}\right.
  \label{eq:gimp-indicator-function}
\end{align}
The spatial gradient of the basis functions is defined as
\begin{align}
  \grad S_{i,mp} := \frac{1}{V_{mp}}\int_{\Omega_{mp}\cap\Omega} \chi_{mp}(\xi) \grad N_{i}(\xi)\,\od\xi\,.
  \label{eq:gimp-basis-gradient}
\end{align}
The product of 1D basis functions gives multi-dimensional shape functions, as in the FEM.
Using these basis functions, we approximate the trial solution fields as
\begin{align}
    \tensor{u}^{h} = \tensor{S}_{u}\cdot\tensor{U}\,, \quad
    {p}^{h} = \tensor{S}_{p}\cdot\tensor{P}\,, \quad
    \label{eq:trial-solution-interpolation}
\end{align}
and their time derivatives as
\begin{align}
    \tensor{v}^{h} = \tensor{S}_{u}\cdot\dot{\tensor{U}}\,, \quad
    \tensor{a}^{h} = \tensor{S}_{u}\cdot\ddot{\tensor{U}}\,,
    \label{eq:trial-solution-time-derivatives-interpolation}
\end{align}
where the superscript $h$ denote spatially discrete versions of continuous fields, $\tensor{S}_{u}$ and $\tensor{S}_{p}$ are arrays of the above-described basis functions interpolating the displacement field and the pore pressure field, respectively,
and $\tensor{U}$ and $\tensor{P}$ are nodal vectors for the displacement field and the pore pressure field, respectively.
We will use linear shape functions for both $\tensor{S}_{u}$ and $\tensor{S}_{p}$.
The variations are approximated by the same basis functions for the trial solutions.

Following the standard discretization procedure, we arrive at the vector form of the variational equations.
At this point, we write them as residual vectors so they can be used to develop a nonlinear solution scheme later on.
For an element $e$, the residual vector of the momentum balance equation is given by
\begin{align}
  [\tensor{\mathcal{R}}_{\rm mom}]^{i}_{e}
  :=& - \int_{\Omega_e} (\symgrad\tensor{S}_{u}^{i})^{\sf T}:\tensor{\sigma}'_{\rm inv}(\tensor{u}^{h})\, \od v
  + \int_{\Omega_e} (\diver\tensor{S}_{u}^{i})^{\sf T}p^{h}\, \od v
  - \int_{\Omega_e}(\symgrad\tensor{S}_{u}^{i})^{\sf T}:\tensor{\sigma}'_{\rm vis}(\tensor{v}^{h})\, \od v \nonumber\\
  &- \int_{\Omega_e} (\tensor{S}_{u}^{i})^{\sf T}\cdot(\rho \tensor{a}^{h})\, \od v
  + \int_{\Omega_e} (\tensor{S}_{u}^{i})^{\sf T}\cdot(\rho\tensor{g})\, \od v
  + \int_{\partial\Omega_e} (\tensor{S}_{u}^{i})^{\sf T}\cdot\hat{\tensor{t}}\,\od a\,,
  \label{eq:residual-momentum}
\end{align}
and the residual vector of the mass balance equation is given by
\begin{align}
  [\tensor{\mathcal{R}}_{\rm mass}]^{i}_{e}
  :=& \int_{\Omega_e} (\tensor{S}_{p}^{i})^{\sf T} (\diver\tensor{v}^{h})\, \od v
  - \int_{\Omega_e} (\grad\tensor{S}_{p}^{i})^{\sf T}\cdot\tensor{q}(\tensor{a}^{h},p^{h})\, \od v
  - \int_{\partial\Omega_e} (\tensor{S}_{p}^{i})^{\sf T}\hat{q}\,\od a\,.
  \label{eq:residual-mass}
\end{align}
For quasi-static problems, the two terms in Eq.~\eqref{eq:residual-momentum} that contain \revised{$\tensor{\sigma}'_{\rm vis}$} and $\tensor{a}^{h}$ vanish, and the second term in Eq.~\eqref{eq:residual-mass} becomes no longer dependent on $\tensor{a}^{h}$. The latter is because $\tensor{a}=\tensor{0}$ in Eq.~\eqref{eq:darcy-law} under quasi-static conditions.

The above equations will be discretized using a single type of material point that is a mixture of the solid and fluid phases.
The main reason is that the interactions between the solid and fluid phases in the porous material are smeared in the poromechanical formulation.
Yet it is worth nothing that there is an alternative approach that uses separate material points for the solid and fluid phases~\cite{Mackenzie-Helnwein2010,Abe2014,Bandara2015}.
This approach can be used for explicit modeling of fluid--solid interactions in porous media; however, its computational cost is significantly larger as it requires a far larger number of material points.
See~\cite{Yerro2015,Soga2016} for review and comparison of these two approaches.

\subsection{Temporal discretization}
We further discretize the semi-discrete residual vectors in time, using different methods for dynamic and quasi-static formulations.
For the dynamic formulation, we use the Newmark family of time integration methods for hyperbolic problems.
Let us denote by $\Del{t} := t_{n+1} - t_{n}$ the time increment from $t_{n}$ to $t_{n+1}$ and write quantities at the new time instance $t_{n+1}$ without subscripts for notational brevity.
The Newmark time integration algorithms can be written as
\begin{align}
  \tensor{v}^{h} &\approx
  \frac{\gamma}{\beta\Del{t}}(\tensor{u}^{h} - \tensor{u}^{h}_{n}) + \left(1 - \frac{\gamma}{\beta}\right)\tensor{v}^{h}_{n}
  + \left(1 - \frac{\gamma}{2\beta}\right)\Del{t}\tensor{a}^{h}_{n}\,,\\
  \tensor{a}^{h} &\approx
  \frac{1}{\beta\Del{t}^{2}}(\tensor{u}^{h} - \tensor{u}^{h}_{n}) - \frac{1}{\beta\Del{t}}\tensor{v}_{n} + \left(1 - \frac{1}{2\beta}\right)\tensor{a}_{n}\,,
  \label{eq:newmark}
\end{align}
where $\beta$ and $\gamma$ are time integration parameters controlling the accuracy and stability of the integration algorithm.
The Newmark algorithms are unconditionally stable for $2\beta\geq\gamma\geq0.5$, and here we choose $\beta=0.3025$ and $\gamma=0.6$.

For the quasi-static formulation, it suffices to discretize $\diver\tensor{v}^{h}$ in time.
In a quasi-static problem where the configuration may significantly change from time $t_{n}$ to $t_{n+1}$, the spatial divergence operator may also considerably evolve during the time step.
For this reason, it would be more accurate to replace $\diver\tensor{v}^{h}$ with $\dot{J}/J$, as done in Choo~\cite{Choo2018d} for an updated Lagrangian formulation for large-deformation poromechanics.
Therefore, using the implicit Euler method which is first-order accurate and unconditionally stable, we replace $\tensor{v}^{h}$ in Eq.~\eqref{eq:residual-mass} with
\begin{align}
  \diver\tensor{v}^{h} = \frac{\dot{J}}{J} \approx \frac{\ln J - \ln J_{n}}{\Del{t}}\,.
  \label{eq:implicit-euler}
\end{align}
Inserting these time-discretized variables into Eqs.~\eqref{eq:residual-momentum} and~\eqref{eq:residual-mass}, we obtain the fully-discrete variational equations.
The discrete equations have two unknown vectors, namely the nodal displacement vector, $\tensor{U}$, and the nodal pore pressure vector, $\tensor{P}$, for both dynamic and quasi-static formulations.

\subsection{Evaluation of the discretized equations}
Because the MPM uses an updated Lagrangian approach, the discretized variational equations are reckoned with respect to the current configuration at time $t_{n+1}$.
However, in the solution step (Stage 2 in Fig.~\ref{fig:mpm-procedure}), the coordinate system is defined with respect to the configuration at the previous time $t_{n}$.
Therefore, for the evaluation of the discretized equations, some quantities have to be pulled back to the previous configuration, which may be thought as the reference configuration in an updated Lagrangian framework.

For this purpose, let $\Omega_{n}$ and $\partial\Omega_{n}$ denote the domain and the boundary of the configuration at time $t_{n}$ calculated from the converged solution in the previous load step.
Likewise, let us denote by $\tensor{F}_{n}$ the deformation gradient at time $t_{n}$ and by $J_{n} := \od V_{n}/\od V$ the Jacobian at time $t_{n}$ which maps the differential volume in the initial configuration, $\od V$, to the differential volume at the previous time instance, $\od V_{n}$.
Then we define the relative deformation gradient and the relative Jacobian as
\begin{align}
    \Del{\tensor{F}} &:= \tensor{F}\cdot\tensor{F}_{n}^{-1} = \tensor{1} + \grad_{n}(\tensor{u} - \tensor{u}_{n})\,, \\
    \Del{J} &:= J/J_{n} = \od v/\od V_{n}\,,
\end{align}
where $\grad_{n}$ is the gradient operator evaluated with respect to the configuration at time $t_{n}$.
Then the differential volume and the material point volume in the current configuration can be calculated as
\begin{align}
    \od v &= \Del{J}\,\od V_{n}\,, \\
    V_{mp} &= \Del{J}(V_{mp})_{n}\,.
\end{align}
Also, we can compute the spatial gradient as
\begin{align}
  \grad(\circ) = \grad_{n}(\circ)\cdot(\Del{\tensor{F}})^{-1}\,,
\end{align}
The symmetric gradient and divergence operators can also be calculated in this way.
Note that all quantities and operators at time $t_{n}$ can be evaluated with respect to a known, fixed coordinate system.
% Also, the differential area can be pulled back using Nanson's formula, given by $\tensor{n}\, \od a = \Del{J}\tensor{N}_{n}\cdot(\Del{\tensor{F}})^{-1}\, \od A_{n}$.

Once the above equations are inserted into the discrete residual vectors, one can readily calculate the discretized governing equations in the MPM.
It is noted that the resulting equations are very similar to the variational equations formulated in~\cite{Borja1995a,Borja1998,Li2004,Uzuoka2012,Borja2016} for mixed finite elements for large-deformation poromechanics.
The key difference is that the equations are now reckoned in the previous configuration, instead of the initial configuration, because the material configuration is evolving in every step in the MPM.
Accordingly, the Jacobian and the deformation gradient in the finite element formulations have been replaced by their relative ones, namely $\Del{J}$ and $\Del{F}$.
Indeed, the relative deformation gradient, $\Del{F}$, has already been commonly used for updating stresses in large-deformation formulations, see Chapter 5 of Borja~\cite{Borja2013} for example.
Therefore, the above approach enables one to develop a mixed MPM formulation for poromechanics as a straightforward extension of mixed finite elements for large-deformation poromechanics.

\subsection{Solution of the coupled system of equations}
Since a large-deformation problem always involves geometric nonlinearity, we use Newton's method to solve the system of equations at hand.
We drive Newton iterations until the residual vector meets the criterion of
\begin{align}
  \frac{\|\tensor{\mathcal{R}}^{k}\|}{\|\tensor{\mathcal{R}}^{0}\|} \leq 10^{-8}\,,
  \label{eq:newton-tolerance}
\end{align}
where $\|\tensor{\mathcal{R}}\|$ is the $L^{2}$-norm of the entire residual vector, with a subscript $k$ denoting the Newton iteration counter.

To calculate an increment vector in each Newton iteration, we construct a Jacobian matrix by linearizing the residual vectors defined in Eqs.~\eqref{eq:residual-momentum} and~\eqref{eq:residual-mass}.
The Jacobian system for the Newton iteration takes the form of
\begin{align}
  \left[\begin{array}{ll}
    \tensor{\mathcal{A}} & \tensor{\mathcal{B}}_{1} \\
    \tensor{\mathcal{B}}_{2} & \tensor{\mathcal{C}}
  \end{array}\right]
  \left\{\begin{array}{ll}
    \Del\tensor{U} \\
    \Del\tensor{P}
  \end{array}\right\}
  =
  - \left\{\begin{array}{l}
    \tensor{\mathcal{R}}_{\rm mom} \\
    \tensor{\mathcal{R}}_{\rm mass}
  \end{array}\right\}\,,
  \label{eq:jacobian-system}
\end{align}
where $\Del\tensor{U}$ and $\Del\tensor{P}$ are Newton increments for nodal vectors of $\tensor{u}^{h}$ and $p^{h}$, and $\tensor{\mathcal{A}}$, $\tensor{\mathcal{B}}_{1}$, $\tensor{\mathcal{B}}_{2}$, and $\tensor{\mathcal{C}}$ are the sub-matrices in the Jacobian matrix, given by
\begin{align}
    \tensor{\mathcal{A}} := \delta_{\tensor{u}^{h}}\tensor{\mathcal{R}}_{\rm mom}\,, \quad
    \tensor{\mathcal{B}}_{1} := \delta_{p^{h}}\tensor{\mathcal{R}}_{\rm mom}\,, \quad
    \tensor{\mathcal{B}}_{2} := \delta_{\tensor{u}^{h}}\tensor{\mathcal{R}}_{\rm mass}\,, \quad
    \tensor{\mathcal{C}} := \delta_{p^{h}}\tensor{\mathcal{R}}_{\rm mass}\,,
\end{align}
where $\delta(\circ)$ is the linearization operator with subscripts denoting the variables linearized.
More specific expressions for these sub-matrices are lengthy but they are essentially the same as those developed for mixed finite element formulations for dynamic and quasi-static poromechanics at large strains~\cite{Li2004}.
For brevity, we omit these expressions and refer to the literature.
Note that a few terms that relate to inertial and damping effects are unwieldy to linearize, but their effects on the overall convergence behavior are not critical as shown in Li \etal~\cite{Li2004}.
We also do not linearize these terms and obtain nearly optimal convergence rates, which will be demonstrated in a numerical example later.
It is also noted that for the dynamic formulation, the time derivatives of $\tensor{u}^{h}$ and $p^{h}$ integrated by the Newmark methods~\eqref{eq:newmark} are linearized as
\begin{align}
    \delta_{\tensor{u}^{h}}\tensor{a}^{h} = \frac{1}{\beta\Del{t}^{2}}\,, \quad
    \delta_{\tensor{u}^{h}}\tensor{v}^{h} = \frac{\gamma}{\beta\Del{t}}\,,
    \label{eq:newmark-linearization}
\end{align}
and for the quasi-static formulation, the Jacobian term is linearized as
\begin{align}
    \delta_{\tensor{u}^{h}}\ln J = \diver\tensor{u}^{h}\,.
\end{align}
Therefore, for both dynamic and quasi-static formulations, we have to solve a $2\times2$ block-partitioned system of linear equations in the form of Eq.~\eqref{eq:jacobian-system} in each Newton iteration.

To solve the coupled linear system~\eqref{eq:jacobian-system} during Newton iterations, we make use of a fully-implicit Krylov method aided by the preconditioning strategy presented in White \etal~\cite{White2016}.
This scheme solves the coupled system in a monolithic way, but applies preconditioners to the Schur complement $\tensor{\mathcal{S}}:=\tensor{\mathcal{C}} - \tensor{\mathcal{B}}_{2}\tensor{\mathcal{A}}^{-1}\tensor{\mathcal{B}}_{1}$ and then $\tensor{\mathcal{A}}$ in a sequential way.
See White \etal~\cite{White2016} for more details and analyses.
We particularly use algebraic multigrid preconditioners for sub-preconditioners and the biconjugate gradient stabilized method for the solver.
This block-preconditioned solution method has performed well for the fully-implicit MPM simulations of various poromechanical problems.

The fully-implicit Newton--Krylov method described above is one of the most robust and efficient solution strategies for coupled poromechanical problems that involve nonlinearity.
However, it is noted that the MPM formulation can be solved by any other method for nonlinear problems or can be time-discretized by an explicit method.
This remains the case even after applying the stabilization method described in the next section.

\subsection{Post-solution update}
Once the above-described solution procedure is finished (\ie~Eq.~\eqref{eq:newton-tolerance} is satisfied), we update state variables of material points (Stage 3 in Fig.~\ref{fig:mpm-procedure}).
For the dynamic formulation, we update $\tensor{u}^{h}$, $p^{h}$, and their first and second time derivatives of the material points, namely $\tensor{v}^{h}$, $\tensor{a}^{h}$, $\dot{p}^{h}$, and $\ddot{p}^{h}$.
Note that, although $\dot{p}^{h}$ does not appear when the pore fluid is incompressible, it will be introduced in stabilization terms later, and its integration via the Newmark algorithm requires $\ddot{p}^{h}$.
For the quasi-static formulation, it suffices to update $\tensor{u}$ and $p$.
All the variables can be updated as  previously described in Eq.~\eqref{eq:map-from-nodes-to-particles},
except for the solid velocity vector, $\tensor{v}^{h}$.
The update scheme in Eq.~\eqref{eq:map-from-nodes-to-particles}, which comes from the PIC method~\cite{Harlow1964}, is known to introduce excessive numerical damping for dynamic problems~\cite{Stomakhin2013,Nairn2015,Hammerquist2017,Jiang2017,Liang2019}.
Therefore,  we use the FLIP method~\cite{Brackbill1986} to update $\tensor{v}^{h}$ at material point $mp$ (denoted by $\tensor{v}_{mp}$) as
\begin{align}
  \tensor{v}_{mp} = (\tensor{v}_{mp})_{n} + \Del{t} \sum_{i=1}^{n_{d}} S_{i,mp} [(1 - \gamma)(\tensor{a}_{i})_{n} + \gamma\tensor{a}_{i}]\,.
\end{align}
Compared with the PIC method, the FLIP method is known to be less stable in time.
For this reason, previous studies have proposed and used more advanced velocity update schemes, such as combined PIC and FLIP methods~\cite{Stomakhin2013,Nairn2015,Liang2019}, the extended PIC method~\cite{Hammerquist2017}, the affine PIC method~\cite{Jiang2017}.
Note, however, that this temporal stability is a completely different issue from the inf--sup stability condition which is the subject of this paper.
Because the stabilization method proposed in this work is independent of the choice of a velocity update scheme,
here we simply use the FLIP method which is most standard in dynamic MPM formulations.

Subsequently, we update the positions of the material points---and also the influence domains if the GIMP method is used---according to the new state variables (Stage 4 in Fig.~\ref{fig:mpm-procedure}).
For updating the lengths of the influence domains, here we follow the method suggested by Charlton \etal~\cite{Charlton2017} whereby the symmetric material stretch tensor is used to prevent spurious disappearance of domains by large rotational deformation.
The background mesh may also be updated at this stage if necessary, although we will simply use a fixed background mesh in this work.

We then move to the next time step and begin the new step by mapping the state variables to nodes (Stage 1 in Fig.~\ref{fig:mpm-procedure}) from the material points at new positions.

% SECTION 4
% ------------------------------------------------------------------------------
\section{Stabilized material point method formulations}
\label{sec:stabilization}
In this section, we present stabilized versions of the MPM formulations for which standard low-order interpolation functions can be used throughout the entire range of drainage conditions.
We begin by explaining the stability condition arising in the undrained limit.
Subsequently, we develop stabilization terms for the dynamic and quasi-static MPM formulations using the PPP method.
The specific forms of the stabilization terms are different for dynamic and quasi-static problems because of time-integration methods and parameters.
Lastly, we describe how the stabilization terms can be implemented particularly for the GIMP method.

\subsection{Stability condition for undrained poromechanical problems}
A coupled poromechanical problem is in an undrained condition when the permeability of the material is sufficiently small for a given time interval.
In this case, there is virtually no relative flow between the pore fluid and the solid deformation (\ie~$\tensor{q}\approx \tensor{0}$), and so the $\tensor{\mathcal{C}}$ block in the Jacobian matrix becomes nearly zero.
Then the coupled linear system~\eqref{eq:jacobian-system} becomes
\begin{align}
  \left[\begin{array}{ll}
    \tensor{\mathcal{A}} & \tensor{\mathcal{B}}_{1} \\
    \tensor{\mathcal{B}}_{2} & \tensor{0}
  \end{array}\right]
  \left\{\begin{array}{ll}
    \Del\tensor{U} \\
    \Del\tensor{P}
  \end{array}\right\}
  =
  - \left\{\begin{array}{l}
    \tensor{\mathcal{R}}_{\rm mom} \\
    \tensor{\mathcal{R}}_{\rm mass}
  \end{array}\right\}\,.
  \label{eq:jacobian-system-undrained}
\end{align}
This matrix structure shows that in the undrained limit, the pore pressure field acts as a Lagrange multiplier imposing the incompressibility constraint for the solid matrix.
The same matrix structure also arises in other types of constrained problems, notably the Stokes flow and incompressible elasticity problems.

For constrained problems, it is well known that the combination of discrete spaces for the primary unknown field and the constraint field should satisfy the inf--sup condition for numerical stability~\cite{Brezzi1974,Brezzi1990}.
In the context of poromechanics, the inf--sup condition can be written as
\begin{align}
  \sup_{\substack{\tensor{\eta}^{h}\in\tensor{\mathcal{U}}^{h}\\ \tensor{\eta}\neq\tensor{0}}}
  \frac{\psi^{h}\diver\tensor{\eta}^{h}}{\|\tensor{\eta}^{h}\|_{1}}
  \geq C\|\psi^{h}\|_{0}\,, \quad
  \forall \psi^{h}\in\mathcal{P}^{h}\,,
  \label{eq:inf--sup}
\end{align}
where $\tensor{\mathcal{U}}^{h}$ and $\mathcal{P}^{h}$ are the discrete spaces for the displacement and pore pressure interpolation functions, respectively,
and $C$ is a positive constant independent of the element size.
When the spaces for interpolation functions satisfy this condition, the discrete constrained problem is well-posed.
Otherwise, however, the discrete problem may be ill-posed and its numerical solution often manifests checkerboard-like spurious oscillations in the pore pressure field.
This numerical problem in the undrained limit of poromechanics has long been studied in the context of mixed finite elements, see, \eg~\cite{Vermeer1981,Murad1994}.

When the displacement and pore pressure fields are interpolated by functions of the same order, the stability condition is not satisfied unfortunately.
Therefore, the MPM cannot use its standard low-order interpolation functions for both the displacement and pore pressure fields, and the use of higher-order interpolation functions for the displacement field is impractical.
In what follows, we develop stabilized formulations that enable one to employ low-order interpolation functions for both fields notwithstanding this stability condition.

\subsection{Stabilized formulations for dynamic and quasi-static poromechanics}
In this work, we use the PPP method to develop stabilized MPM formulations.
This method was originally developed for mixed finite elements for the Stokes flow problem by Bochev \etal~\cite{Bochev2006a} and then applied to mixed finite elements for other constrained problems including undrained poromechanics~\cite{White2008,Sun2013a,Choo2015,Sun2015,Monforte2019,Choo2019}.
The PPP method draws on that equal-order linear interpolation functions satisfy a weaker inf--sup condition, given by
\begin{align}
  \sup_{\substack{\tensor{\eta}\in\tensor{\mathcal{U}}^{h}\\ \tensor{\eta}\neq\tensor{0}}}
  \frac{\psi^{h}\diver\tensor{\eta}^{h}}{\|\tensor{\eta}\|_{1}}
  \geq C_{1}\|\psi^{h}\|_{0} - C_{2}\|\psi^{h}-\Pi \psi^{h}\|_{0}\,, \quad
  \forall p^{h}\in\mathcal{P}^{h}\,,
  \label{eq:weak-inf--sup}
\end{align}
with $C_{1}>0$ and $C_{2}>0$ independent of $h$, and $\Pi:L_{2}(\Omega)\mapsto R_{0}$ a projection operator from the $L^{2}$ space to the piecewise constant space.
Comparing Eqs.~\eqref{eq:inf--sup} and~\eqref{eq:weak-inf--sup}, one can identify that the last term in Eq.~\eqref{eq:weak-inf--sup}, $C_{2}\|\psi^{h}-\Pi \psi^{h}\|_{0}$, is the deficiency of the equal-order linear interpolation functions for inf--sup stability.
To compensate this deficiency, the PPP method augments a stabilization term to the mass balance equation such that
\begin{align}
  \left[\begin{array}{ll}
    \tensor{\mathcal{A}} & \tensor{\mathcal{B}}_{1} \\
    \tensor{\mathcal{B}}_{2} & \tensor{\mathcal{C}} + \tensor{\mathcal{C}}_{\rm stab}
  \end{array}\right]
  \left\{\begin{array}{ll}
    \Del\tensor{U} \\
    \Del\tensor{P}
  \end{array}\right\}
  =
  - \left\{\begin{array}{l}
    \tensor{\mathcal{R}}_{\rm mom} \\
    \tensor{\mathcal{R}}_{\rm mass} + \tensor{\mathcal{R}}_{\rm stab}
  \end{array}\right\}\,.
  \label{eq:jacobian-system-stabilized-ppp}
\end{align}
Here, $\tensor{\mathcal{R}}_{\rm stab}$ is the stabilization term and its linearization gives rise to additional terms $\tensor{\mathcal{C}}_{\rm stab} := \delta_{p^{h}}\tensor{\mathcal{R}}_{\rm stab}$ in the lower diagonal block of the Jacobian matrix.
This term makes the lower diagonal block non-zero even when $\tensor{\mathcal{C}} \approx \tensor{0}$ in the undrained limit.

Specifically, the PPP stabilization term can be developed as follows.
We first define the projection operator $\Pi$ as
\begin{align}
  \Pi \psi^{h}|_{\Omega_e} = \frac{1}{V_{e}}\int_{\Omega_e} \psi^{h}\,\od v\,, \quad
  \label{eq:projection-operator}
\end{align}
with $V_{e}$ the volume of element $e$.
In words, the projection operator calculates the element-wise averages of a discrete field as piecewise constants.
Then the element-wise stabilization term can be written as
\begin{align}
  [\tensor{\mathcal{R}}_{\rm stab}]_{e} := \tau \int_{\Omega_{e}}(\psi^{h} - \Pi\psi^{h})(\dot{p}^{h} - \Pi\dot{p}^{h}) \,\od v,
  \label{eq:stab-term-ppp}
\end{align}
where $\tau$ is a stabilization parameter which will be discussed below in more detail.
Here, we have written the stabilization term as a semi-discrete form because the time derivative of $\dot{p}^{h}$ will be discretized differently for dynamic and quasi-static problems.

For dynamic problems, we use the Newmark algorithm so replace $\dot{p}^{h}$ in Eq.~\eqref{eq:stab-term-ppp} with
\begin{align}
  \dot{p}^{h} &=
  \frac{\gamma}{\beta\Del{t}}(p^{h} - p^{h}_{n}) + \left(1 - \frac{\gamma}{\beta}\right)\dot{p}^{h}_{n}
  + \left(1 - \frac{\gamma}{2\beta}\right)\Del{t}\ddot{p}^{h}_{n}\,, \\
  \ddot{p}^{h} &=
  \frac{1}{\beta\Del{t}^{2}}(p^{h} - p^{h}_{n}) - \frac{1}{\beta\Del{t}}\dot{p}^{h}_{n} + \left(1 - \frac{1}{2\beta}\right)\ddot{p}^{h}_{n}\,.
  \label{eq:p-dot-newmark}
\end{align}
Note that, whereas $\ddot{p}^{h}$ does not appear in the formulation, it is necessary to evaluate $\dot{p}^{h}$ using the Newmark algorithm.
For this reason, material points should also carry $\dot{p}^{h}$ and $\ddot{p}^{h}$ as their state variables, as described in the previous section.
However, because both $\dot{p}^{h}$ and $\ddot{p}^{h}$ can be linearized with respect to $p^{h}$ (\cf~Eq.~\eqref{eq:newmark-linearization} for linearization of $\tensor{v}^{h}$ and $\tensor{a}^{h}$), the stabilization term does not introduce any additional degree of freedom.

For quasi-static problems, we use the implicit Euler method and discretize $\dot{p}^{h}$ in Eq.~\eqref{eq:stab-term-ppp} as
\begin{align}
  \dot{p}^{h} = \frac{p^{h} - p^{h}_{n}}{\Del{t}}\,.
  \label{eq:p-dot-implicit-euler}
\end{align}
This term can also be linearized with respect to $p^{h}$.

\revised{
\subsection{Stabilization parameter}
The stabilization parameter, $\tau$, in Eq.~\eqref{eq:stab-term-ppp} determines the degree of stabilization effects.
If this parameter is too small, the numerical solution may still suffer from the inf--sup stability problem and manifest checkerboard oscillations.
Conversely, too large a stabilization parameter may render the numerical solution overly diffusive.
Therefore, it is critical to use an appropriately large stabilization parameter that successfully suppresses oscillations without undesirable smoothing.

For mixed FEM for coupled poromechanics, a few different expressions have been proposed for $\tau$.
When White and Borja~\cite{White2008} first applied the PPP method to poromechanical problems under infinitesimal strain and quasi-static conditions, they used
\begin{align}
  \tau = \tau_{\rm W} := \frac{1}{2G}\,,
  \label{eq:tau-white}
\end{align}
where $G$ the shear modulus of the solid matrix.
This expression is obtained by replacing the fluid viscosity in the original PPP method for the Stokes flow problem~\cite{Bochev2006a} with the shear modulus of the solid matrix.
Later on, Sun \etal~\cite{Sun2013a} extended the PPP method to large-deformation quasi-static poromechanics and derived an expression for $\tau$ from a simple 1D setting.
The expression is given by
\begin{align}
  \tau = \tau_{\rm S} := \frac{1}{M}\left(1 - 3\frac{c_v\Del{t}}{h^2}\right)\left[1 + \tanh\left(2 - 12\frac{c_v\Del{t}}{h^2} \right)\right]\,,
  \label{eq:tau-sun}
\end{align}
where $M = K + (4/3)G$ is the P-wave modulus, $h$ is the element size, and $c_{v} = M\kappa$ is the coefficient of consolidation.
More recently, for dynamic poromechanics, Monforte \etal~\cite{Monforte2019} proposed the following expression
\begin{align}
  \tau = \tau_{\rm M} := \max\left(\frac{2\kappa}{c_v} - \frac{\beta}{\gamma}\frac{12\kappa\Del{t}}{h^2},0\right)\,.
  \label{eq:tau-monforte}
\end{align}
Here, the maximum operator is used to prevent a negative stabilization parameter.
Notably, this expression uses the Newmark algorithm parameters, $\beta$ and $\gamma$, so it is only suitable for dynamic problems.

It is noted that all of the above three expressions can be calculated from material properties and time integration parameters, without any tuning for a specific problem.
Also, at least for the problems we tested, the three expressions lead to similar stabilization performance for both dynamic and quasi-static problems.
This finding will be demonstrated later in Section~\ref{sec:examples}.
Therefore, although the effect of the stabilization term is sensitive to the value of $\tau$,
any of the existing expressions for $\tau$ may be used at one's discretion.
}

\subsection{Implementation of stabilization terms}
For both dynamic and quasi-static formulations, the stabilization terms can be implemented in an existing MPM code without change in the existing spatial discretization scheme or the time integration algorithm.
This is a significant advantage over the fractional step method used in Jassim \etal~\cite{Jassim2013} which requires a substantial change in the time stepping method.
This advantage may be the main reason that non-operator-splitting stabilization methods have recently been developed for mixed finite elements for dynamic poromechanical problems~\cite{dePouplana2017,Monforte2019}, let alone quasi-static problems.

The only non-trivial issue arises when the stabilization term is implemented in conjunction with the GIMP method.
This is because the stabilization term~\eqref{eq:stab-term-ppp} contains element-wise averages of the shape functions and the pore pressure solutions.
For the standard MPM, these values can be evaluated in the same way as in the FEM.
However, in the GIMP method, a material point has an influence domain that may affect multiple elements, so the average within an element in the background mesh is not so meaningful.
To address this issue, here we adopt an approach suggested by Coombs \etal~\cite{Coombs2018} for the implementation of the $\bar{\tensor{F}}$ method for the GIMP method.
In essence, it introduces a basis function that is constant in an influence domain, given by
\begin{align}
  S_{i,mp}^{0} := \frac{1}{V_{mp}}\int_{\Omega_{mp}} \frac{1}{2}\,\od\xi\,.
\end{align}
Then, for the average (projected) terms in Eq.~\eqref{eq:stab-term-ppp}, this basis function is used in lieu of the original GIMP basis function defined in Eq.~\eqref{eq:gimp-basis}.
See Coombs \etal~\cite{Coombs2018} for more details.
Note that this still does not introduce significant additional efforts for implementation.
Importantly, this approach allows us to calculate the stabilization term even when an element has a single or few material points, which is a crucial feature under very large deformations.

% SECTION 5
% ------------------------------------------------------------------------------
\section{Numerical examples}
\label{sec:examples}
This section presents four numerical examples that verify, investigate, and demonstrate the performance of the stabilized MPM formulations.
The purpose of the first and second examples is to verify and investigate the stabilization formulations using quasi-static and dynamic problems previously simulated by the FEM.
The third and fourth examples are intended to demonstrate the capability of the stabilized MPM formulations for simulating problems that are challenging for mesh-based Lagrangian methods.

In each example, we will compare results obtained by the standard (non-stabilized) MPM and stabilized MPM formulations.
For undrained poromechanical problems, the standard MPM formulation is expected to show spurious oscillations in the pore pressure field,
whereas the stabilized formulations are developed to suppress these oscillations.
Therefore, in what follows, we focus on the comparison of pore pressure fields obtained by these formulations.
We will particularly consider excess pore pressure, which is the transient portion of pore pressure.
\revised{Also, unless otherwise specified, we use $\tau_{\rm W}$ in Eq.~\eqref{eq:tau-white} for the stabilization parameter, $\tau$.}

\revised{
The proposed stabilization methods are applicable to both the original MPM and the GIMP method.
In the first example where particles do not cross any cell, the original and GIMP methods give the same numerical solution.
However, when the original MPM was used for other problems, we encountered cell-crossing effects that are irrelevant to oscillations due to inf--sup stability.
To demonstrate the coexistence of inf--sup and cell-crossing instabilities, in the second example we purposely use the original MPM for generating unstable results in the undrained limit.
For the third and fourth examples, we use the GIMP method to obtain both unstable and stable results to fully focus on the inf--sup stability problem.
}

All the MPM formulations have been implemented in \texttt{Geocentric-MPM}, an in-house MPM code for geomechanics built upon the \verb|deal.II| finite element library~\cite{dealII90}, \verb|p4est| mesh handling library~\cite{Burstedde2011}, and the \verb|Trilinos| project~\cite{Heroux2012}.
The code is an MPM extension of a massively parallel finite element code \texttt{Geocentric} used in the literature~\cite{White2011a,Choo2016,White2016,Choo2018b,Choo2018c}.

\subsection{Terzaghi's 1D consolidation problem}
\begin{figure}[htbp]
  \centering
  \includegraphics[width=0.18\textwidth]{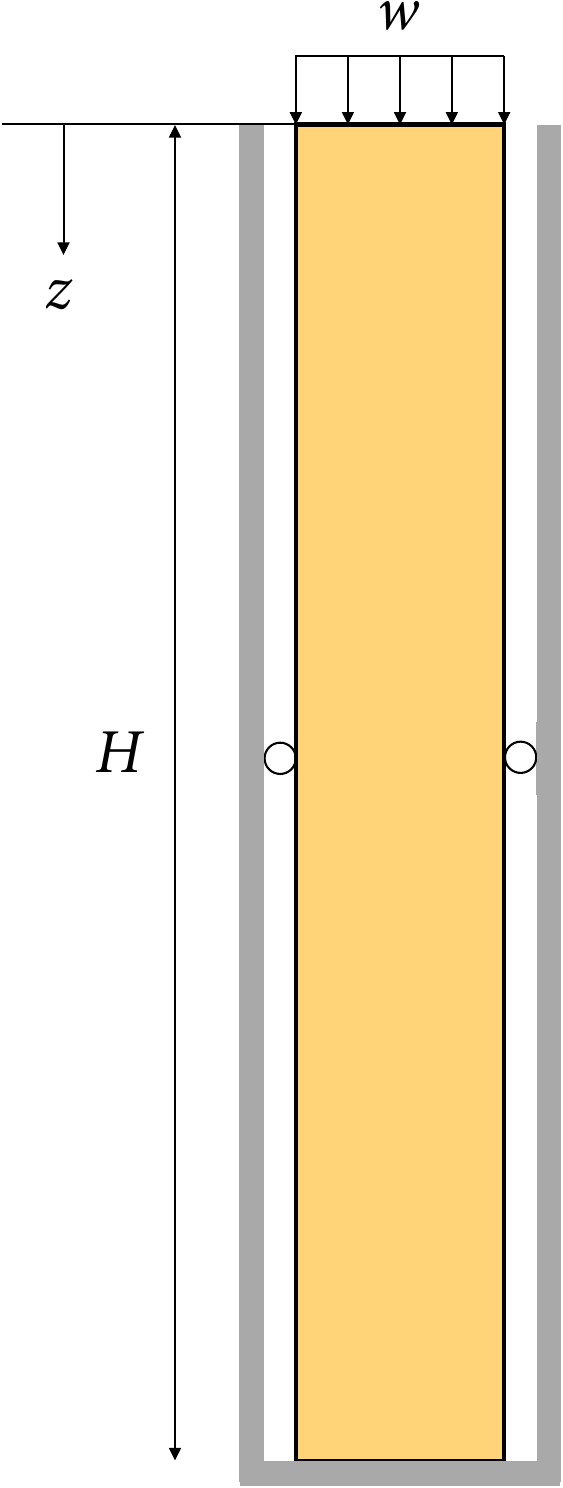}
  \caption{Setup of Terzaghi's 1D consolidation problem.}
  \label{fig:terzaghi-setup}
\end{figure}

Our first example is Terzaghi's 1D consolidation problem, which has commonly been used for the verification and assessment of stabilized mixed finite elements for poromechanics (\eg~\cite{White2008,Preisig2011,Yoon2018,Choo2019}).
The geometry and boundary conditions of this problem are depicted in Fig.~\ref{fig:terzaghi-setup}.
Here, a uniformly distributed static load $w$ is applied on the top boundary of a saturated porous column of height $H$.
The top boundary is drained while all other boundaries are impermeable.
The two lateral boundaries are supported by rollers and the bottom boundary is fixed.
Gravity is neglected in this problem.

This problem permits an analytical solution in the limit of infinitesimal deformation.
To express the analytical solution, we define dimensionless pressure, depth, and time, denoted by $P$, $Z$, and $T$, respectively, as
\begin{align}
  P := p/w\,, \quad Z := z/H\,, \quad T := (c_v/H^2)t\,.
\end{align}
The analytical solution is given by
\begin{align}
  P(Z,T) = \sum_{m=0}^{\infty} \frac{2}{M}\sin(MZ)e^{-M^2T}\,, \quad
  M := \pi(2m+1)/2\,.
  \label{eq:terzaghi-exact}
\end{align}

Using the standard and stabilized MPM formulations under quasi-static conditions, we simulate Terzaghi's problem under an undrained condition.
To get specific solutions, we assign the material's bulk modulus as $K=1$ MPa, Poisson's ratio as $\nu=0.25$, the coefficient of consolidation as $c_{v} = 1.8\times10^{-5}$ m$^2$/s, and the domain height as $H=1$ m.
To compare the numerical solution with the small-strain analytical solution, we apply a small load of $w=0.001$ MPa.
To investigate sensitivity with respect to spatial discretization, we use two levels of discretization: (a) 40 material points with 80 quadrilateral elements in the background mesh, and (b) 80 material points within 160 quadrilateral elements.
Each element contains 2 material points in both cases.
\revised{For the 1D consolidation problem, Vermeer and Verruijt~\cite{Vermeer1981} derived that a time increment smaller than $(1/6)h^{2}/c_{v}$ will lead to unstable numerical solution.
Because our purpose is to examine the performance of the stabilization method, we deliberately set $\Del{t}=0.1$ second which is smaller than $(1/6)h^{2}/c_{v}$ in both cases.
Therefore, unless stabilized, the numerical solution should exhibit checkerboard oscillations.
Then we run a single time step to $t$ = 0.1 second.}

\begin{figure}[htbp]
  \centering
  \subfloat[40 material points\newline]{\includegraphics[width=0.85\textwidth]{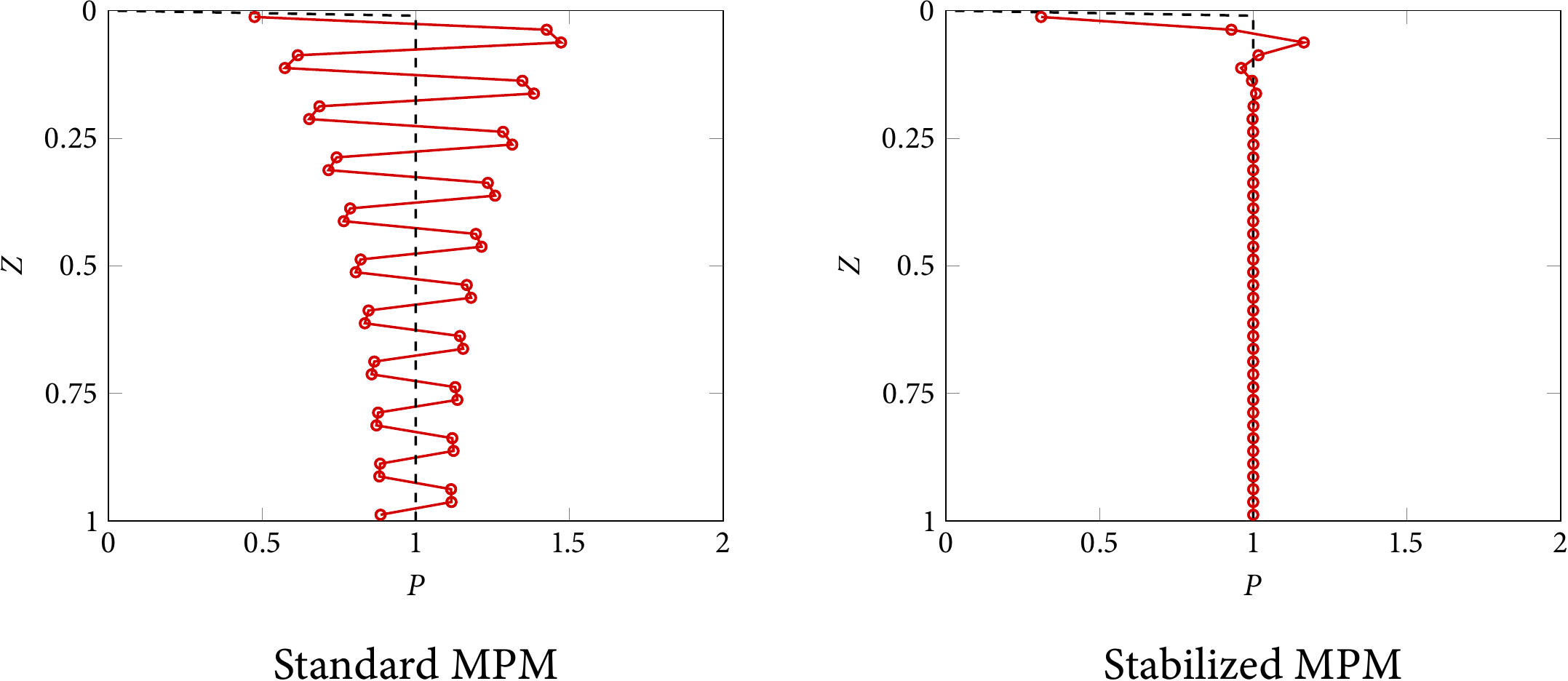}\label{fig:terzaghi-medium}}\\ \vspace{1em}
  \subfloat[80 material points\newline]{\includegraphics[width=0.85\textwidth]{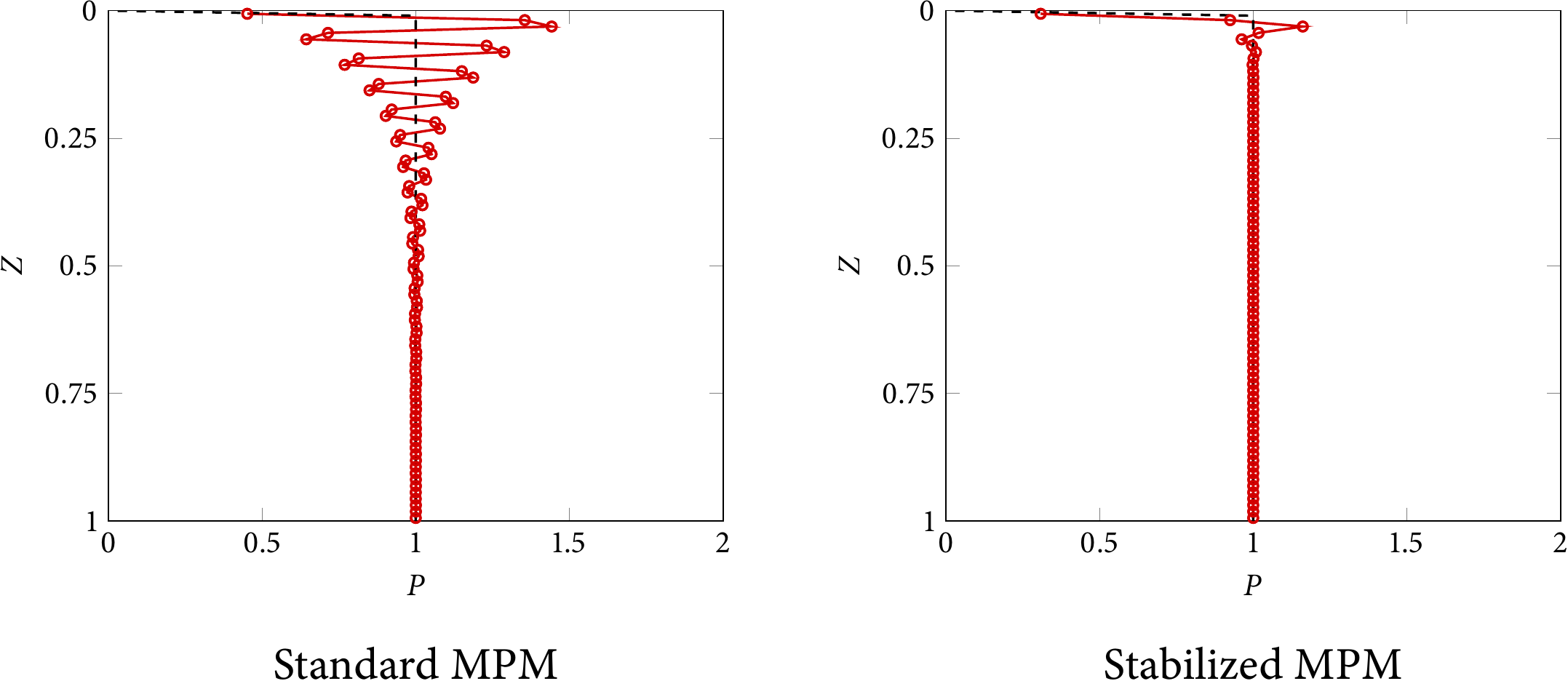}\label{fig:terzaghi-fine}}
  \caption{Excess pore pressure distributions at $t$ = 0.1 second. Dashed lines denote analytical solutions.}
  \label{fig:terzaghi-results}
\end{figure}

Figure~\ref{fig:terzaghi-results} presents the numerical solutions obtained by the standard and stabilized MPM formulations.
The analytical solution is also drawn as dashed lines for reference.
The numerical results of the two MPM formulations are found to be very similar to their FEM counterparts in the literature~\cite{White2008,Choo2019}.
In Fig.~\ref{fig:terzaghi-medium}, the standard MPM formulation without stabilization shows severe pressure oscillation throughout the domain.
The pressure oscillation is somewhat alleviated in a finer discretization (Fig.~\ref{fig:terzaghi-fine}), but it is still unacceptable in the upper half of the domain.
Conversely, the results of the stabilized MPM formulation are free of these oscillations owing to the inf--sup stability problem.
Note that, while the stabilized result still shows some oscillation near the drainage boundary, it is due to an inherent drawback of a continuous Galerkin interpolation rather than the inf--sup stability.
See, for example, Fig. 5 of White and Borja~\cite{White2008} in which intrinsically stable Taylor--Hood elements also show the same oscillations near the drainage boundary.
Because this is also consistent with results in the stabilized FEM literature~\cite{White2008,Choo2019},
the stabilized MPM formulation has been verified for this quasi-static problem.

\begin{figure}[htbp]
  \centering
  \includegraphics[width=0.5\textwidth]{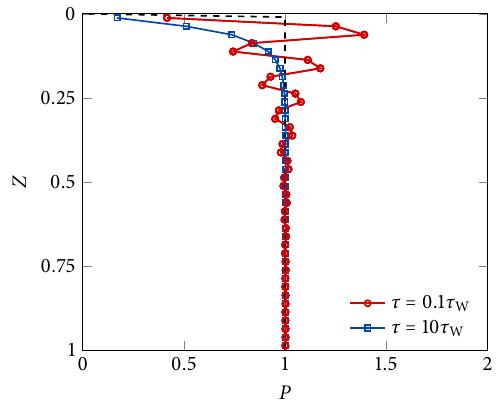}
  \caption{\revised{Comparison of stabilized MPM solutions obtained by two stabilization parameters, $0.1\tau_{\rm W}$ and $10\tau_{\rm W}$, for the 40 material points case.}}
  \label{fig:terzaghi-parameter}
\end{figure}

\begin{figure}[htbp]
  \centering
  \subfloat[40 material points\newline]{\includegraphics[width=0.45\textwidth]{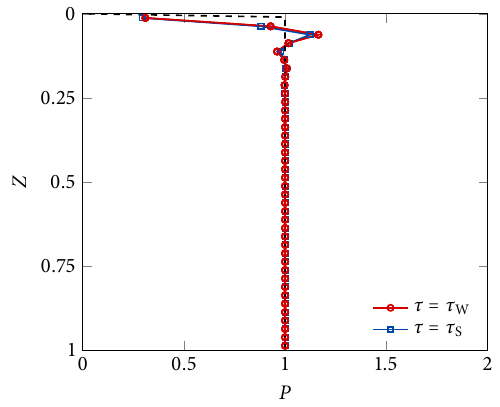}} \hspace{1em}
  \subfloat[80 material points\newline]{\includegraphics[width=0.45\textwidth]{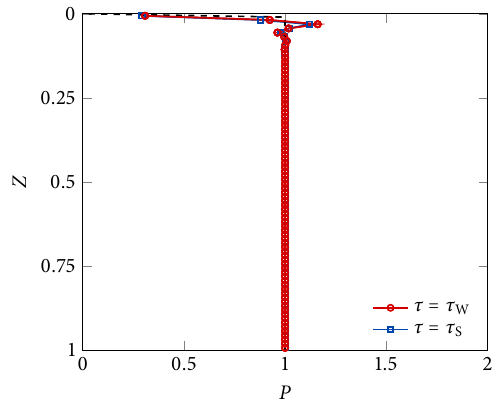}}
  \caption{\revised{Comparison of stabilized MPM solutions obtained by two stabilization parameters, $\tau_{\rm W}$ and $\tau_{\rm S}$.}}
  \label{fig:terzaghi-white-sun}
\end{figure}

\revised{
To investigate the parameter sensitivity of the stabilization method, we repeat the 40 material points case with 0.1 and 10 times of the original stabilization parameter, $\tau_{\rm W}$, and present the results in Fig.~\ref{fig:terzaghi-parameter}.
As can be seen, when the stabilization parameter is not large enough ($\tau=0.1\tau_{\rm W}$), the stabilization effect is insufficient and the checkerboard oscillations persist.
However, if the stabilization parameter is too large ($\tau=10\tau_{\rm W}$), the stabilization term produces excessive artificial diffusion, resulting in over-smoothing in the numerical solution.
This comparison demonstrates the importance of the use of an appropriate stabilization parameter.
We then compare the stabilization performances of two existing expressions for the stabilization parameter for quasi-static problems, $\tau_{\rm W}$ in Eq.~\eqref{eq:tau-white} and $\tau_{\rm S}$ in Eq.~\eqref{eq:tau-sun}, in Fig.~\ref{fig:terzaghi-white-sun}.
It is found that the two stabilization parameters lead to more or less the same numerical solutions.
Notably, although $\tau_{\rm S}$ is a function of the element size, $h$, it does not lead to significant differences from $\tau_{\rm W}$ (which is independent of $h$) for the two discretization cases.
While not presented for brevity, other quasi-static problems that we tested have shown similar insensitivity to the choice between $\tau_{\rm W}$ and $\tau_{\rm S}$.
}


\subsection{Strip footing under harmonic loading}
\begin{figure}[htbp]
  \centering
  \includegraphics[width=0.8\textwidth]{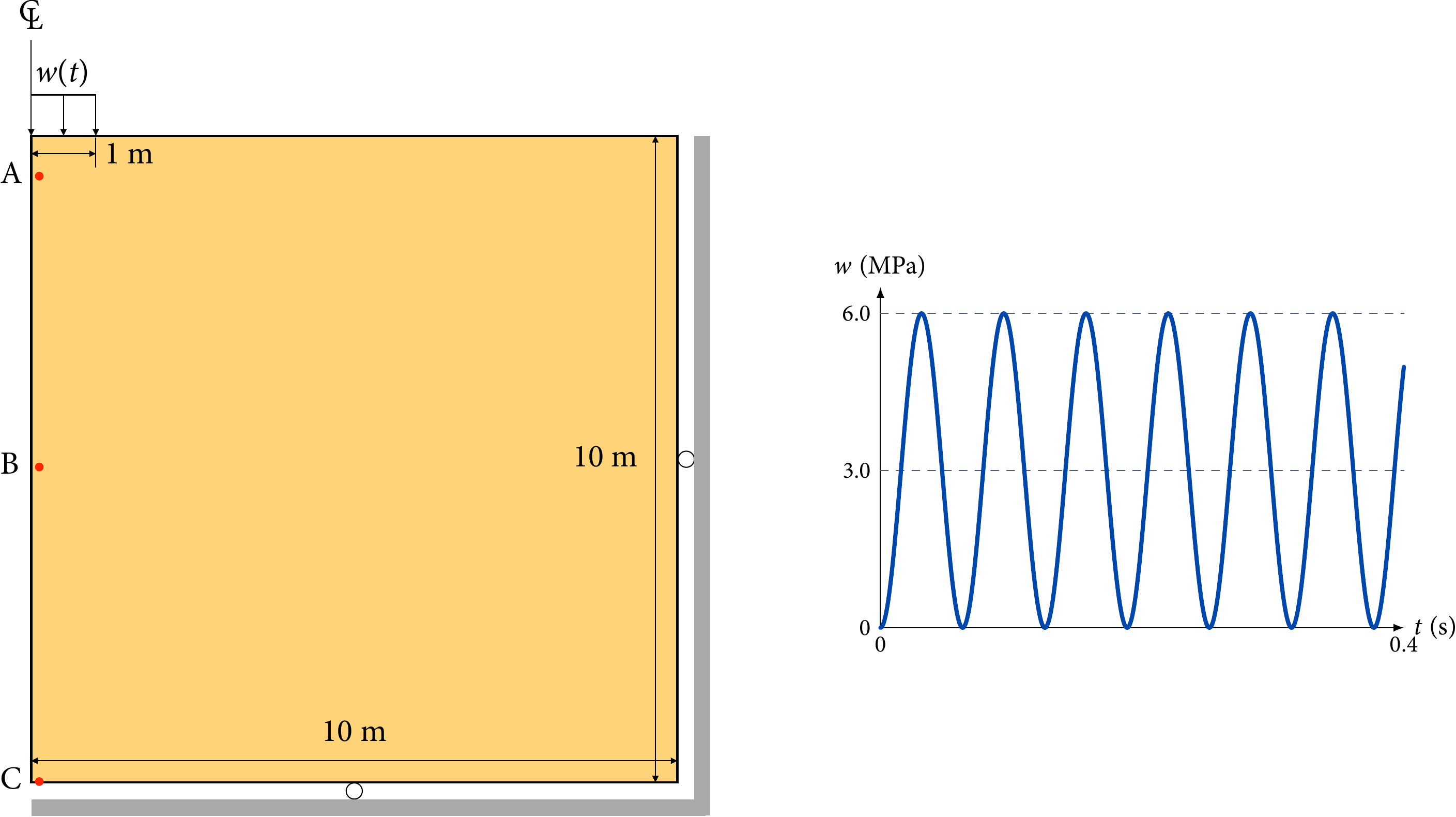}
  \caption{Setup of the strip footing under harmonic loading problem. The loading is $w = 3 - 3\cos(100t)$ MPa ($t$ in seconds). Points A, B, and C are locations at which time evolutions of pore pressures will be presented in Fig.~\ref{fig:harmonic-2d}. The elevations of Points A, B, and C from the bottom boundary are 9.4375 m, 4.9375 m, and 0.0625 m, respectively.}
  \label{fig:harmonic-2d-setup}
\end{figure}

Our next example is a saturated porous ground that underpins a vertically vibrating strip footing.
Figure~\ref{fig:harmonic-2d-setup} illustrates the geometry, boundary conditions, and the loading profile of this strip footing problem.
Most conditions of this problem are similar to those of a numerical example solved by mixed finite elements in Li \etal~\cite{Li2004}.
There are two differences from the original problem: one is that the permeability of the ground has been lowered to a constant value of $10^{-14}$ m$^2$ to introduce undrained deformations, and the other is that a Hencky elasticity model is used instead of a Neo-Hookean model.
Yet the Lam\'{e} parameters, namely $\lambda=8.4$ MPa and $G=5.6$ MPa, are the same as those of the Neo-Hookean model in Li \etal~\cite{Li2004}.
Other parameters are as follows: the damping parameter is $\alpha_{\rm vis} = 0.04$, the dynamic viscosity is $\mu_{f}=10^{-6}$ kPa$\cdot$s, the mass density of the solid phase is $\rho_{s}=2.5$ Mg/m$^3$, the mass density of the fluid phase is $\rho_{w}=1.0$ Mg/m$^3$, and the initial porosity of the material is 0.33.

The domain is 20 m wide and 10 m deep, and its top center is subjected to a 2 m wide strip footing under harmonic loading of $w=3 - 3\cos(100t)$ MPa ($t$ in seconds).
This Neumann boundary condition has been implemented using a method proposed by Bisht and Salgado~\cite{Bisht2018}.
The footing is assumed to be an impermeable boundary while the rest of the top boundary is drained. The bottom boundary and two lateral boundaries are also impermeable.
Taking advantage of symmetry, we model the right half of the domain.

For MPM simulation, we discretize the domain by 1,600 square elements (size $h = 0.25$ m) and fill 4 material points per element.
In addition, to accommodate material points that would heave above the ground surface by the loading, we augment 160 empty elements of the same size onto the top boundary.
We simulate the problem until $t=0.4$ second with a uniform time increment of $\Del{t}=0.001$ second.
\revised{While the size of $\Del{t}$ is not restricted by stability, it is chosen to be small enough to capture the dynamic process at hand with sufficient accuracy.}
Also, for verification of the stabilized MPM results, we simulate the same problem with intrinsically stable Taylor--Hood mixed finite elements~\cite{Taylor1973}.
The Taylor--Hood elements employ quadratic shape functions for $\tensor{u}$ and linear shape functions for $p$, which is a well-known pair that satisfies the inf--sup stability condition in the undrained limit.
For a better comparison of the MPM and FEM solution to this problem, we use Eq.~\eqref{eq:map-from-particles-to-nodes-consistent} for the mapping from material points to nodes, because the MPM's standard mapping method~\eqref{eq:map-from-particles-to-nodes} is known to introduce some numerical damping for a similar cyclic loading problem~\cite{Burgess1992}.
Lastly, for producing standard MPM results without stabilization, we use the original MPM to demonstrate that inf--sup and cell-crossing instabilities coexist when particles are crossing cells.

\begin{figure}[htbp]
  \centering
  \subfloat[$t=0.02$ s\newline]{\includegraphics[width=0.85\textwidth]{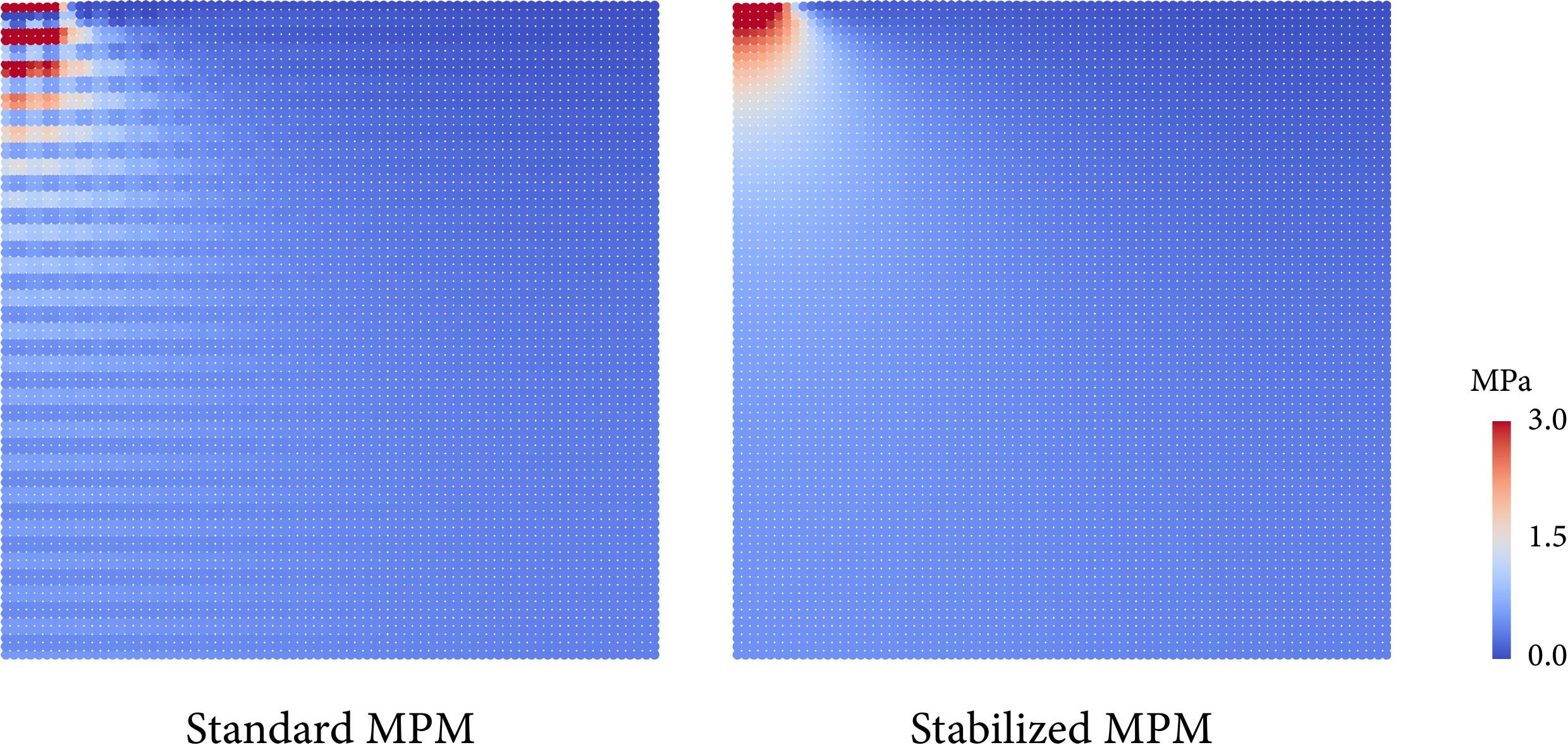}\label{fig:harmonic-2d-contour-t0.02s}}\\ \vspace{1em}
  \subfloat[$t=0.1$ s\newline]{\includegraphics[width=0.85\textwidth]{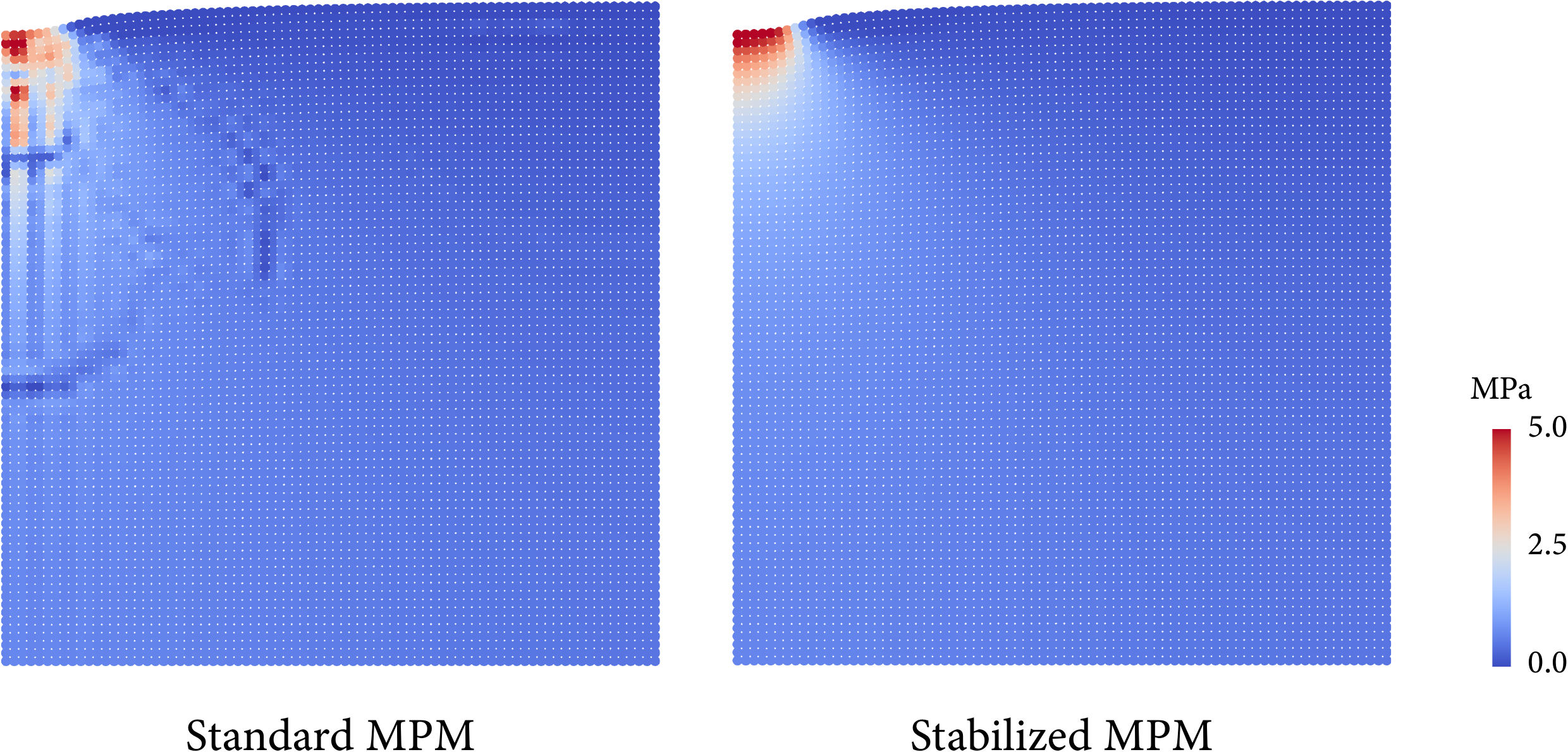}\label{fig:harmonic-2d-contour-t0.1s}}
  \caption{Excess pore pressure fields under harmonic loading.}
\end{figure}

Figure~\ref{fig:harmonic-2d-contour-t0.02s} presents numerical solutions of excess pore pressure fields obtained by the standard and stabilized MPM formulations, after 0.02 second of loading.
We can see that the non-stabilized standard MPM shows checkerboard oscillations in the pressure field, which is a typical manifestation of lack of the inf--sup stability.
By contrast, the stabilized MPM result is free of such oscillation, which confirms the performance of the stabilization method in the undrained limit.
Figure~\ref{fig:harmonic-2d-contour-t0.1s} shows excess pore pressure solutions at a later time, $t=0.1$ second.
Now that the domain has experienced finite deformation, the standard MPM shows two types of oscillations, one due to the inf--sup condition and another due to the cell-crossing instability.
Specifically, the former gives rise to jump-like pressure oscillations right below the footing, and the latter gives rise to line-like oscillations away from the footing.
Regardless of the coexistence of two different instabilities, however, we have found that the combination of the PPP and GIMP methods can provide oscillation-free solutions.

\begin{figure}[htbp]
  \centering
  \subfloat[Point A\newline]{\includegraphics[width=0.75\textwidth]{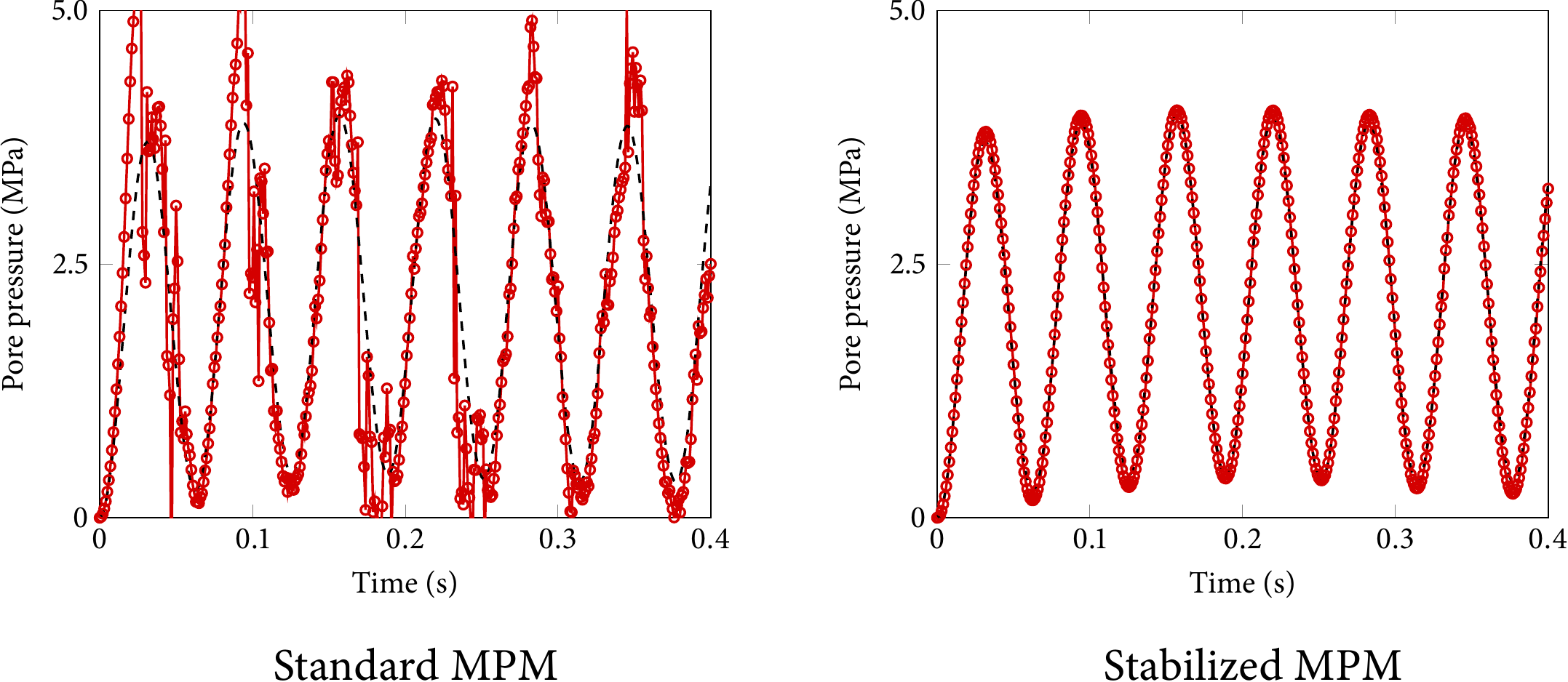}\label{fig:harmonic-2d-top}}\\ \vspace{1em}
  \subfloat[Point B\newline]{\includegraphics[width=0.75\textwidth]{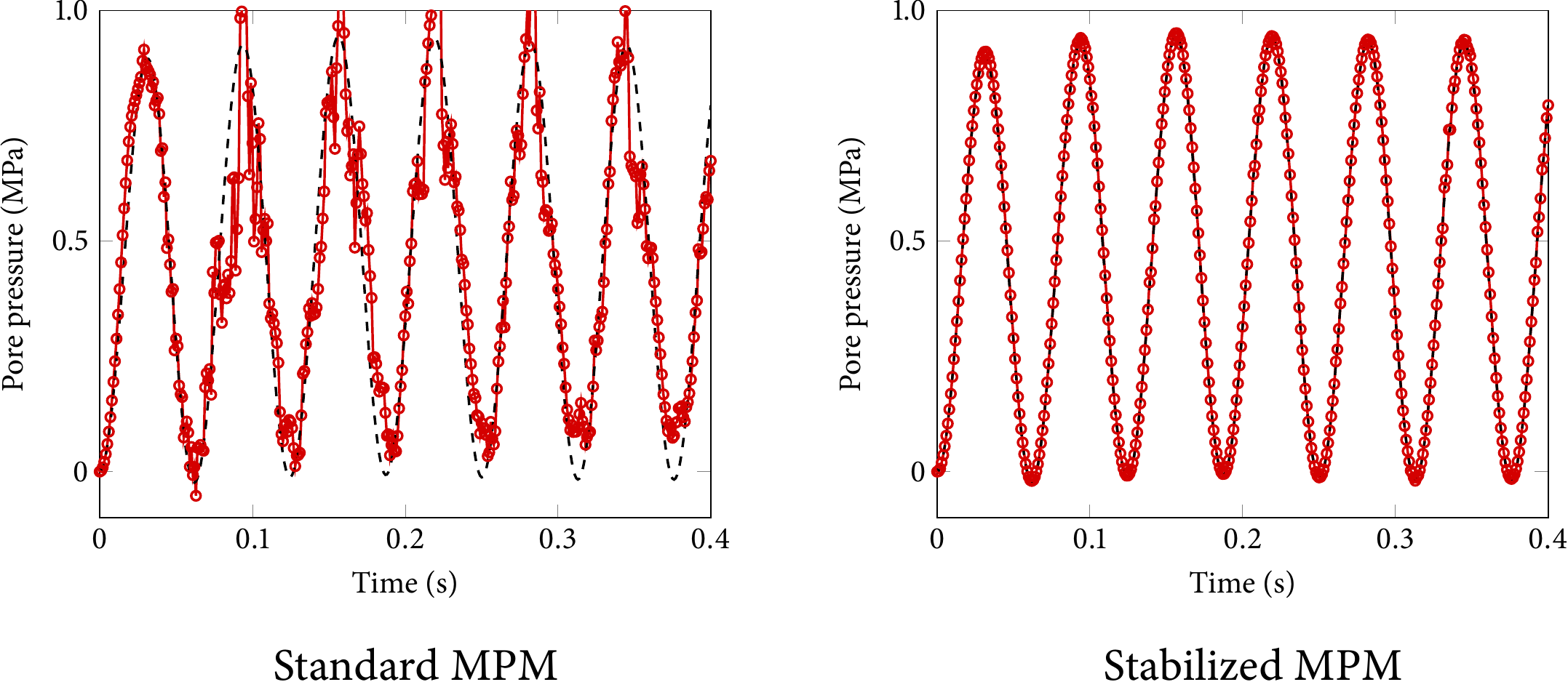}\label{fig:harmonic-2d-middle}}\\ \vspace{1em}
  \subfloat[Point C\newline]{\includegraphics[width=0.75\textwidth]{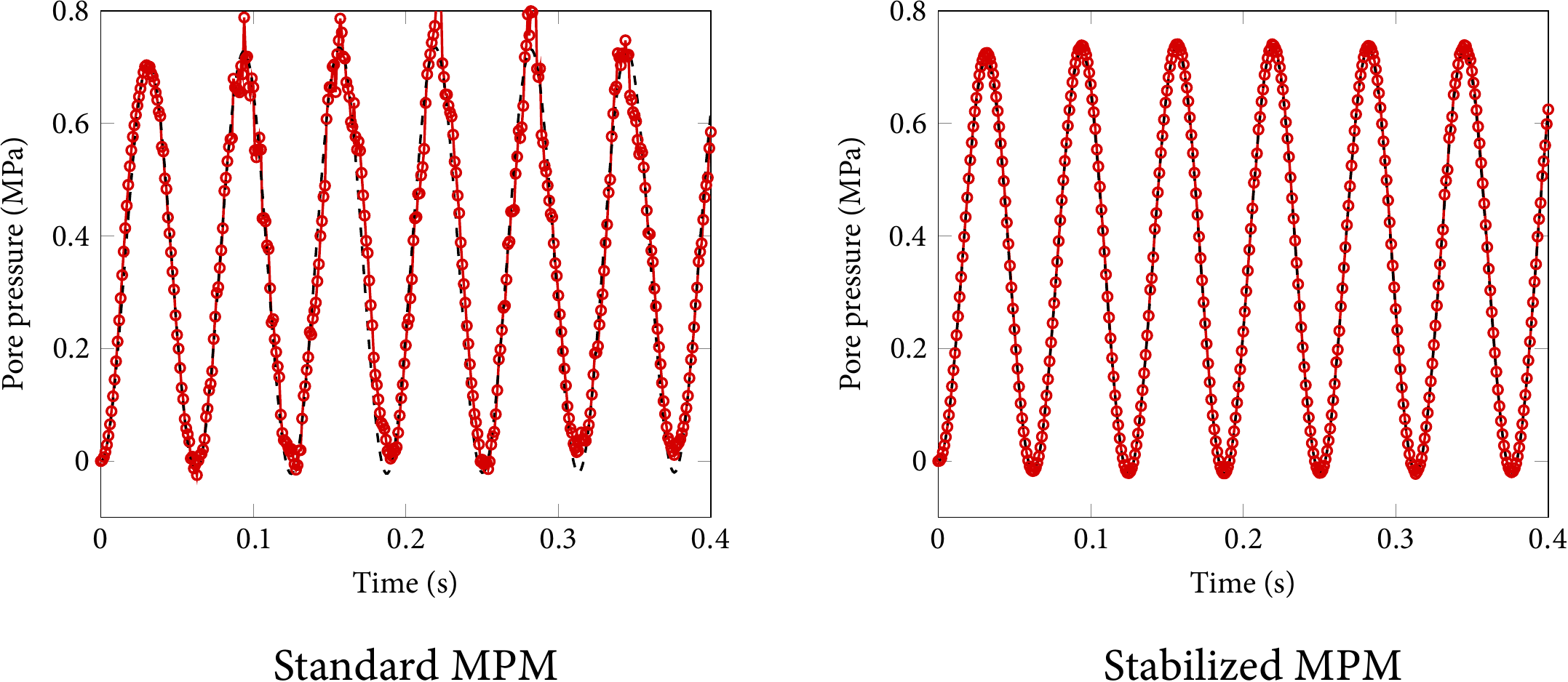}\label{fig:harmonic-2d-bottom}}
  \caption{Time evolution of excess pore pressures at the three points designated in Fig.~\ref{fig:harmonic-2d-setup}. Dashed lines denote numerical solutions obtained by the FEM using stable (Taylor--Hood) elements.}
  \label{fig:harmonic-2d}
\end{figure}

To verify the results of the stabilized MPM formulation, in Fig.~\ref{fig:harmonic-2d} we compare time evolutions of excess pore pressures at three points depicted in Fig.~\ref{fig:harmonic-2d-setup} with solutions obtained by the intrinsically stable Taylor--Hood mixed finite elements.
The results of the standard MPM are also shown for completeness.
It can be seen that the standard MPM gives highly erroneous results, particularly at Point A.
The stabilized MPM, however, allows us to obtain results very close to the stable finite element solutions.

\begin{figure}[htbp]
  \centering
  \subfloat[Point A\newline]{\includegraphics[width=0.4\textwidth]{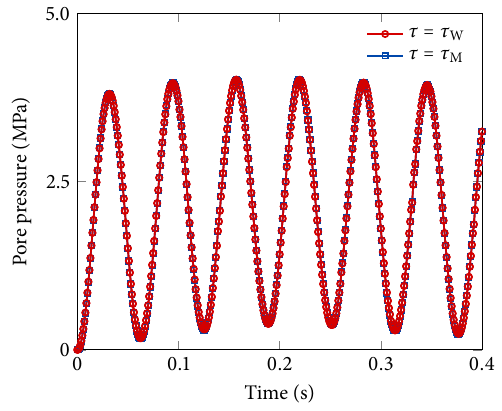}\label{fig:harmonic-2d-top-parameter}}\\ \vspace{1em}
  \subfloat[Point B\newline]{\includegraphics[width=0.4\textwidth]{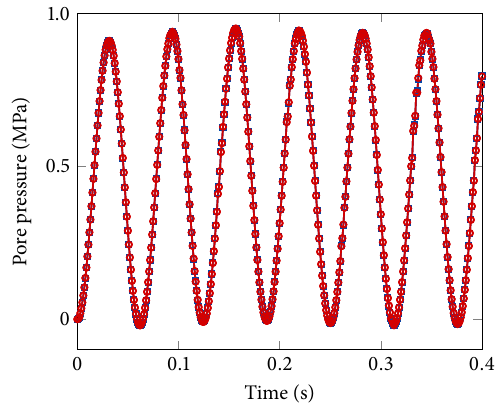}\label{fig:harmonic-2d-middle-parameter}}\\ \vspace{1em}
  \subfloat[Point C\newline]{\includegraphics[width=0.4\textwidth]{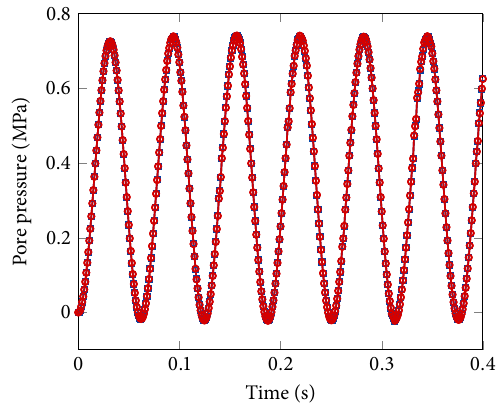}\label{fig:harmonic-2d-bottom-parameter}}
  \caption{\revised{Comparison of stabilized MPM solutions obtained by two stabilization parameters, $\tau_{\rm W}$ and $\tau_{\rm M}$, in terms of the time evolution of excess pore pressures at the three points designated in Fig.~\ref{fig:harmonic-2d-setup}.}}
  \label{fig:harmonic-2d-parameter}
\end{figure}

\revised{Having verified the stabilized MPM formulation, we also examine its sensitivity to the choice of expressions for $\tau$. We repeat the same problem with $\tau_{\rm M}$ in Eq.~\eqref{eq:tau-monforte} which has been proposed by Monforte \etal~\cite{Monforte2019} for dynamic poromechanical problems.
Figure~\ref{fig:harmonic-2d-parameter} compares the results obtained by $\tau_{\rm W}$ and $\tau_{\rm M}$.
It is found that the two results are indeed indistinguishable except few points.
This finding suggests that the original stabilization parameter proposed by White and Borja~\cite{White2008}, $\tau_{\rm W}$, may still be well used for dynamic problems.
For the record, we also observed that too large a $\tau$ (\eg~10$\tau_{\rm W}$) leads to over-smoothing and decreases in the peak pressure values in Fig.~\ref{fig:harmonic-2d-parameter}.}

Lastly, using this example, we demonstrate the performance of the fully-implicit MPM formulation for coupled poromechanics.
For the stabilized MPM simulations, Table~\ref{tab:harmonic-2d-newton-ppp} presents relative residual norms during Newton iterations and the number of iterations during preconditioned Krylov solution of the linear system~\eqref{eq:jacobian-system} at four time steps.
The table shows that the residual norms has converged quickly within 3 Newton iterations, with nearly optimal rates of convergence.
Such convergence rates have been attained even without a couple of complicated terms related to damping and acceleration, as observed in Li \etal~\cite{Li2004} for finite elements.
Moreover, the number of Krylov iterations are fairly small.
These rapid convergences in Newton and Krylov iterations clearly demonstrate that the current fully-implicit MPM is computationally efficient.
We have observed more or less similar performance for all other examples in this section, so we will omit these details for other problems for brevity.
The combination of this solution efficiency and the unconditional stability of the time integration algorithm renders the fully-implicit MPM one of the most powerful approaches to large-deformation poromechanical problems.
\begin{table}[htbp]
  \centering
  \renewcommand*{\arraystretch}{1.3}
  \begin{tabular}[1.3]{c|cccccccc}
    \hline\hline
    % & \multicolumn{8}{c}{MPM-PPP}\\
    % \cline{2-9}
     & \multicolumn{2}{c}{Step 100} & \multicolumn{2}{c}{Step 200} & \multicolumn{2}{c}{Step 300} & \multicolumn{2}{c}{Step 400}\\
    Iteration & $\|\mathcal{\tensor{R}}^{k}\|/\|\mathcal{\tensor{R}}^{0}\|$ & Krylov & $\|\mathcal{\tensor{R}}^{k}\|/\|\mathcal{\tensor{R}}^{0}\|$ & Krylov & $\|\mathcal{\tensor{R}}^{k}\|/\|\mathcal{\tensor{R}}^{0}\|$ & Krylov & $\|\mathcal{\tensor{R}}^{k}\|/\|\mathcal{\tensor{R}}^{0}\|$ & Krylov\\
    \hline
    0 & \multicolumn{1}{r}{1.0e-00} & 12 & \multicolumn{1}{r}{1.0e-00} & 15 & \multicolumn{1}{r}{1.0e-00} & 13 & \multicolumn{1}{r}{1.0e-00} & 13\\
    1 & \multicolumn{1}{r}{1.2e-03} & 15 & \multicolumn{1}{r}{7.3e-04} & 17 & \multicolumn{1}{r}{2.3e-04} & 15 & \multicolumn{1}{r}{1.6e-03} & 13\\
    2 & \multicolumn{1}{r}{2.3e-06} & 14 & \multicolumn{1}{r}{7.1e-07} & 15 & \multicolumn{1}{r}{1.1e-07} & 14 & \multicolumn{1}{r}{4.1e-06} & 14\\
    3 & \multicolumn{1}{r}{3.4e-09} & -  & \multicolumn{1}{r}{5.4e-10} & -  & \multicolumn{1}{r}{4.7e-11} & -  & \multicolumn{1}{r}{5.6e-09} & -\\
    \hline\hline
  \end{tabular}
  \caption{Performance of the Newton--Krylov solver for the stabilized MPM formulation. The tolerance for Krylov iterations was $10^{-8}\|\mathcal{\tensor{R}}^{k}\|$.}
  \label{tab:harmonic-2d-newton-ppp}
\end{table}

\subsection{Self-weight consolidation}
Our third example has two purposes: (a) to show the performance of the stabilized MPM formulation for a 2D quasi-static problem, and (b) to demonstrate the effectiveness of the unconditionally-stable implicit MPM formulation for modeling of a long-term poromechanical process.
For these purposes, we simulate consolidation of a soft porous material due to its self weight, which has relevance to the settlement of a very soft soil such as a dredged marine deposit.
Given the long time scale of the consolidation process, the quasi-static formulation is an appropriate choice for this example.

\begin{figure}[htbp]
  \centering
  \includegraphics[width=0.4\textwidth]{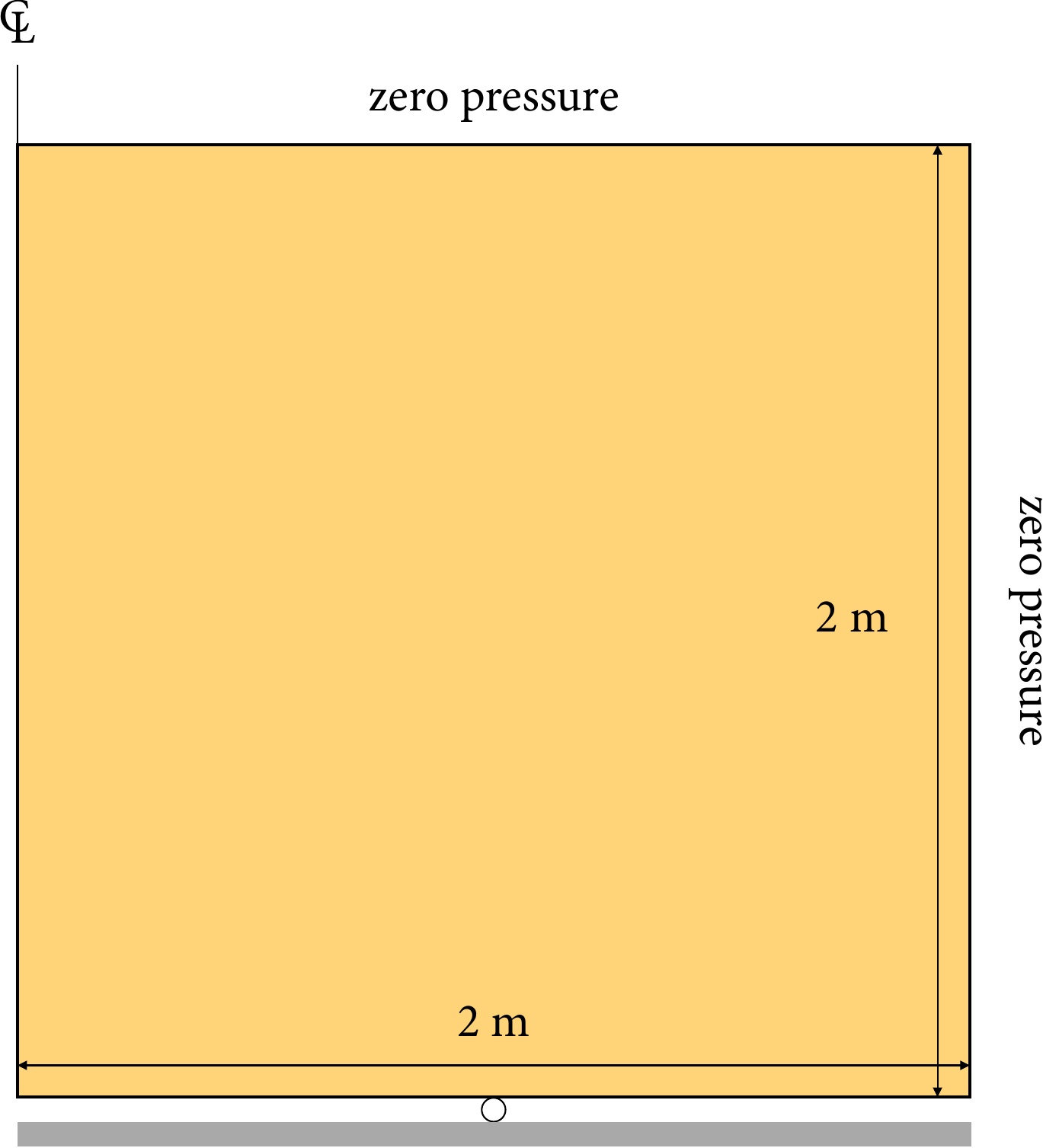}
  \caption{Setup of the self-weight consolidation problem.}
  \label{fig:self-weight-consolidation-setup}
\end{figure}

Figure~\ref{fig:self-weight-consolidation-setup} depicts the setup of this problem.
We consider a 4 m wide and 2 m tall domain, modeling only its right half by taking advantage of symmetry.
No external load is applied because the consolidation process will be driven by the gravitational force applied in the beginning of the problem.
This gravitational force will generate excess pore pressure, and the dissipation of this excess pore pressure will give rise to deformation over time.
The top and right boundaries are zero pressure boundaries through which the pore fluid will be drained during the problem.
The bottom boundary is impermeable and rigid.

As for the material parameters, we assign a bulk modulus of $K = 15$ kPa, a Poisson's ratio of $\nu=0.3$, an initial permeability of 10$^{-14}$ m$^2$, and an initial porosity of 0.5.
In this problem, the permeability is related to porosity change as Eq.~\eqref{eq:kozeny-carman}.
The dynamic viscosity of the fluid is $\mu_{f}=10^{-6}$ kPa$\cdot$s.
The mass densities of the solid and fluid phases are assigned as $\rho_{s}=2.6$ Mg/m$^3$ and $\rho_{w}=1.0$ Mg/m$^3$, respectively.

We begin the problem by applying the gravitational force at $t=0$ and proceed to the consolidation stage with an initial time increment $\Del{t}=0.1$ s.
The time increment is multiplied by a factor of 1.2 after each time step.
\revised{This time stepping scheme is intended to efficiently simulate the increasingly slow consolidation process while still capturing the early phase of the consolidation process accurately.}
To examine the sensitivity of numerical solution to spatial discretization, we consider two levels of discretization: one using 6,400 material points placed on 256 square elements ($h=0.125$ m) in the background grid, and the other using 16,384 material points with 1,024 square elements ($h=0.0625$ m).
For both cases, we simulate the problem until the excess pore pressure becomes fully dissipated.

\begin{figure}[htbp]
  \centering
  \subfloat[6,400 material points\newline]{\includegraphics[width=0.8\textwidth]{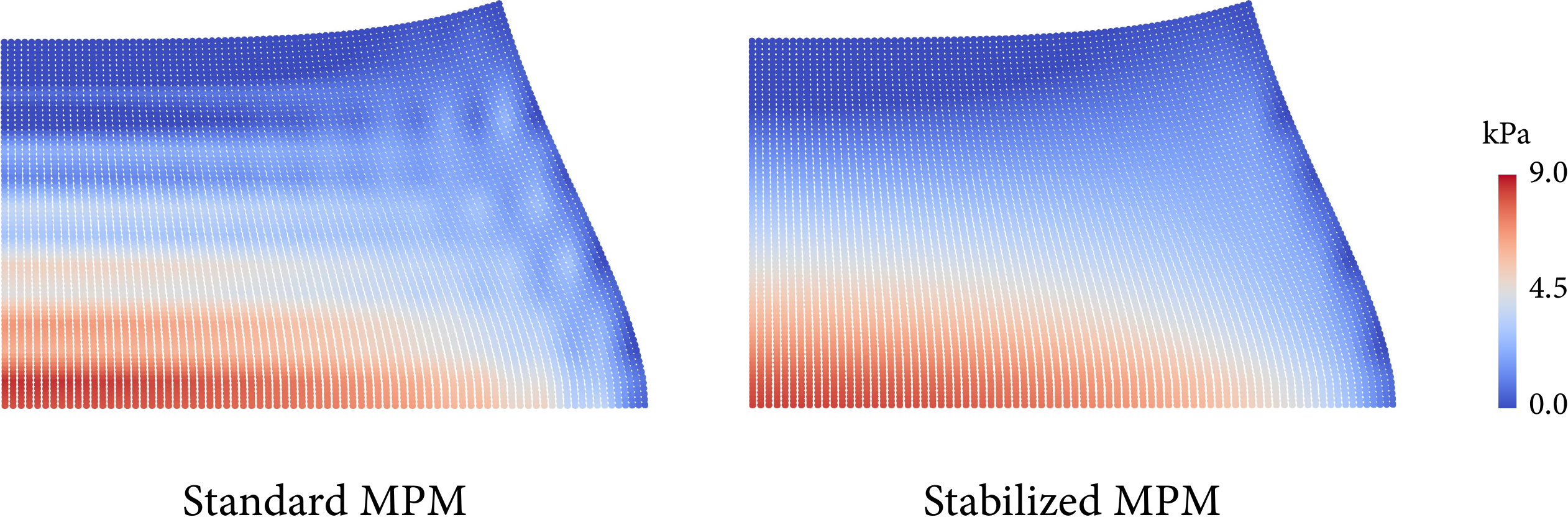}}\vspace{1em}\\
  \subfloat[16,384 material points\newline]{\includegraphics[width=0.8\textwidth]{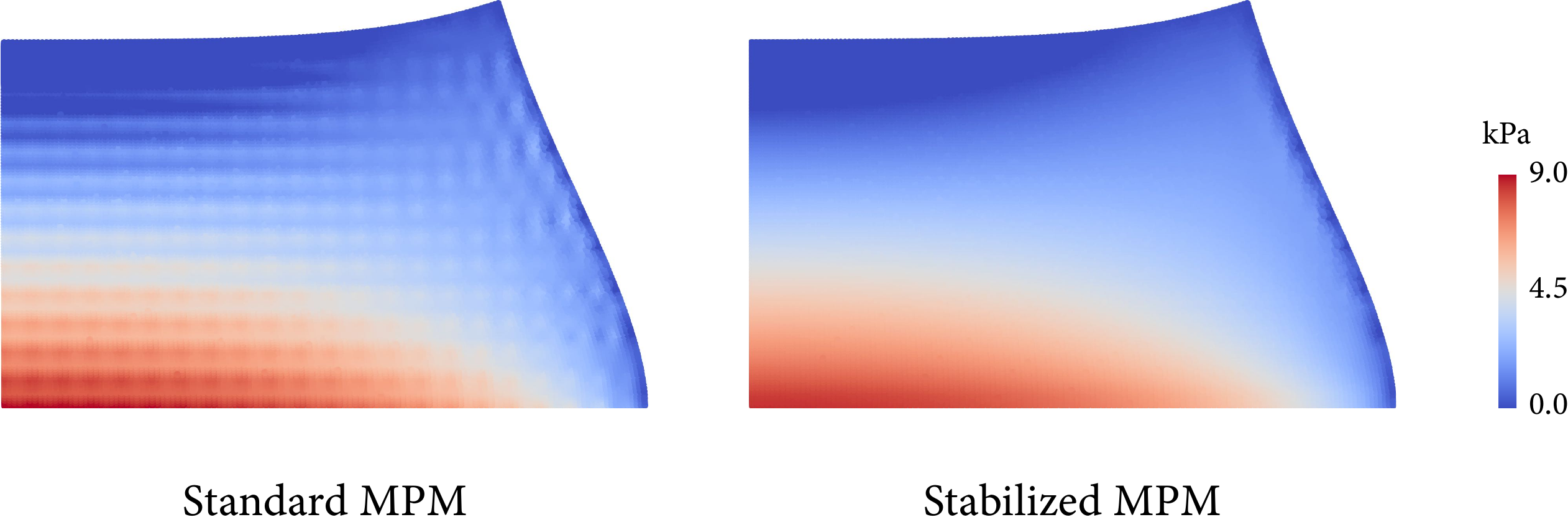}}
  \caption{Excess pore pressure fields immediately after gravity loading.}
  \label{fig:self-weight-consolidation-undrained}
\end{figure}

Figure~\ref{fig:self-weight-consolidation-undrained} shows standard and stabilized MPM solutions of excess pore pressure fields immediately after gravity loading.
One can see that the standard MPM solutions again exhibit checkerboard pressure oscillations, irrespective of the level of spatial discretization.
The stabilized MPM solutions, by contrast, are free of these spurious oscillations and show converging results upon refinement.
It is noted that these solutions have been obtained under a completely undrained condition in which the lower diagonal block of the Jacobian matrix is strictly zero.
Therefore, these results confirm that the stabilized MPM formulation performs well even under the most stringent undrained condition.

\begin{figure}[htbp]
  \centering
  \subfloat[6,400 material points\newline]{\includegraphics[width=0.95\textwidth]{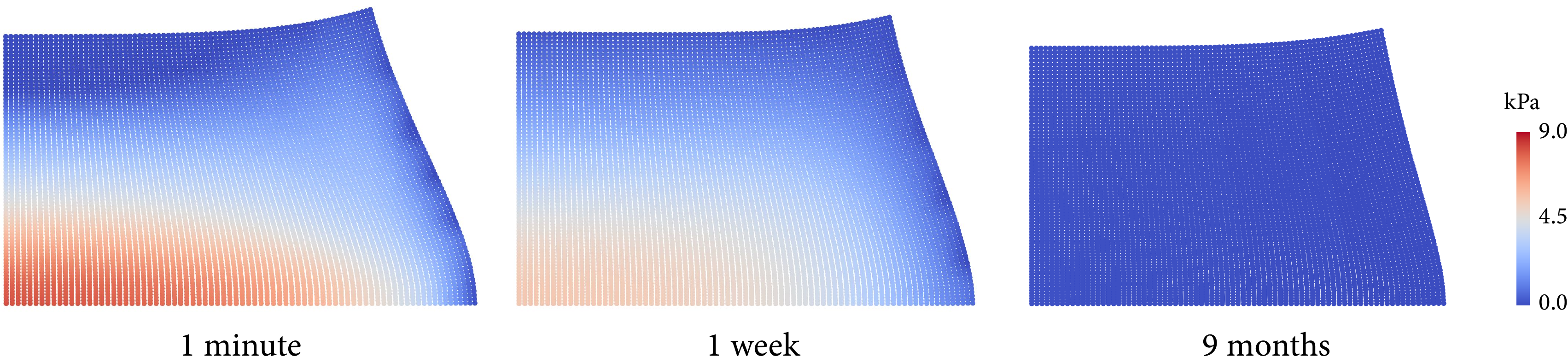}}\vspace{1em}\\
  \subfloat[16,384 material points\newline]{\includegraphics[width=0.95\textwidth]{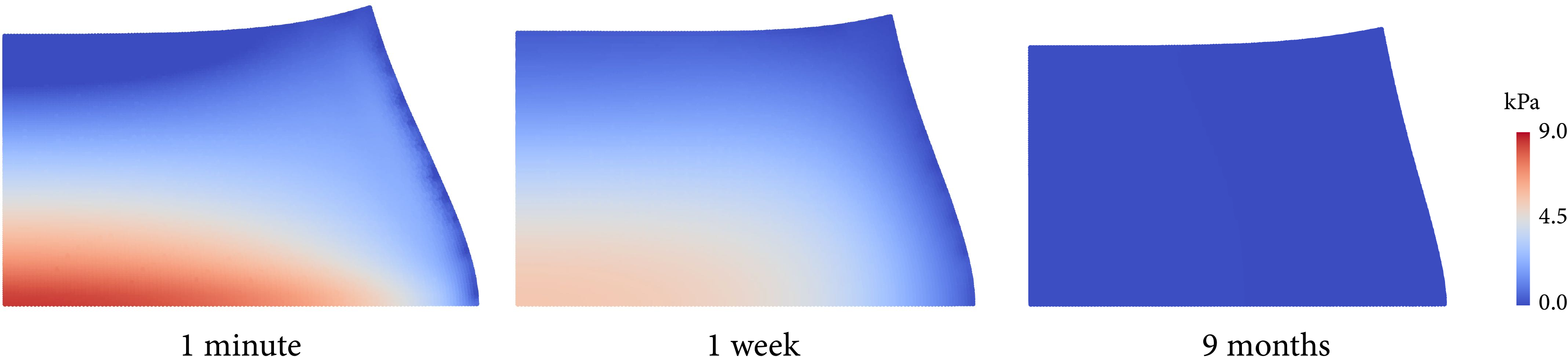}}
  \caption{Time evolution of excess pore pressure fields after gravity loading. All results were obtained by the stabilized MPM.}
  \label{fig:self-weight-consolidation-contours}
\end{figure}

In Fig.~\ref{fig:self-weight-consolidation-contours}, we present simulation results of the consolidation process from the initial conditions in the previous figure.
All these results have been obtained by the stabilized MPM formulation.
In the previous figure, the overall material deformed in an incompressible manner due to the undrained condition.
However, in the current figure, the material shrinks as the pore fluid is discharged into drainage boundaries during the consolidation stage.
We can see that the stabilized MPM also well simulates this consolidation process and associated finite deformations, and that the results again converge upon spatial refinement.
It is worth noting that MPM simulation of this long-term poromechanical process has been made feasible thanks to the use of an unconditionally stable, implicit method.
Because many poromechanical problems in subsurface systems involve long time scales, we believe that the implicit MPM formulation presented in this work will be useful for addressing similar long-term deformation problems.

\subsection{Impact of two poroelastic bodies}
\begin{figure}[htbp]
  \centering
  \includegraphics[width=0.5\textwidth]{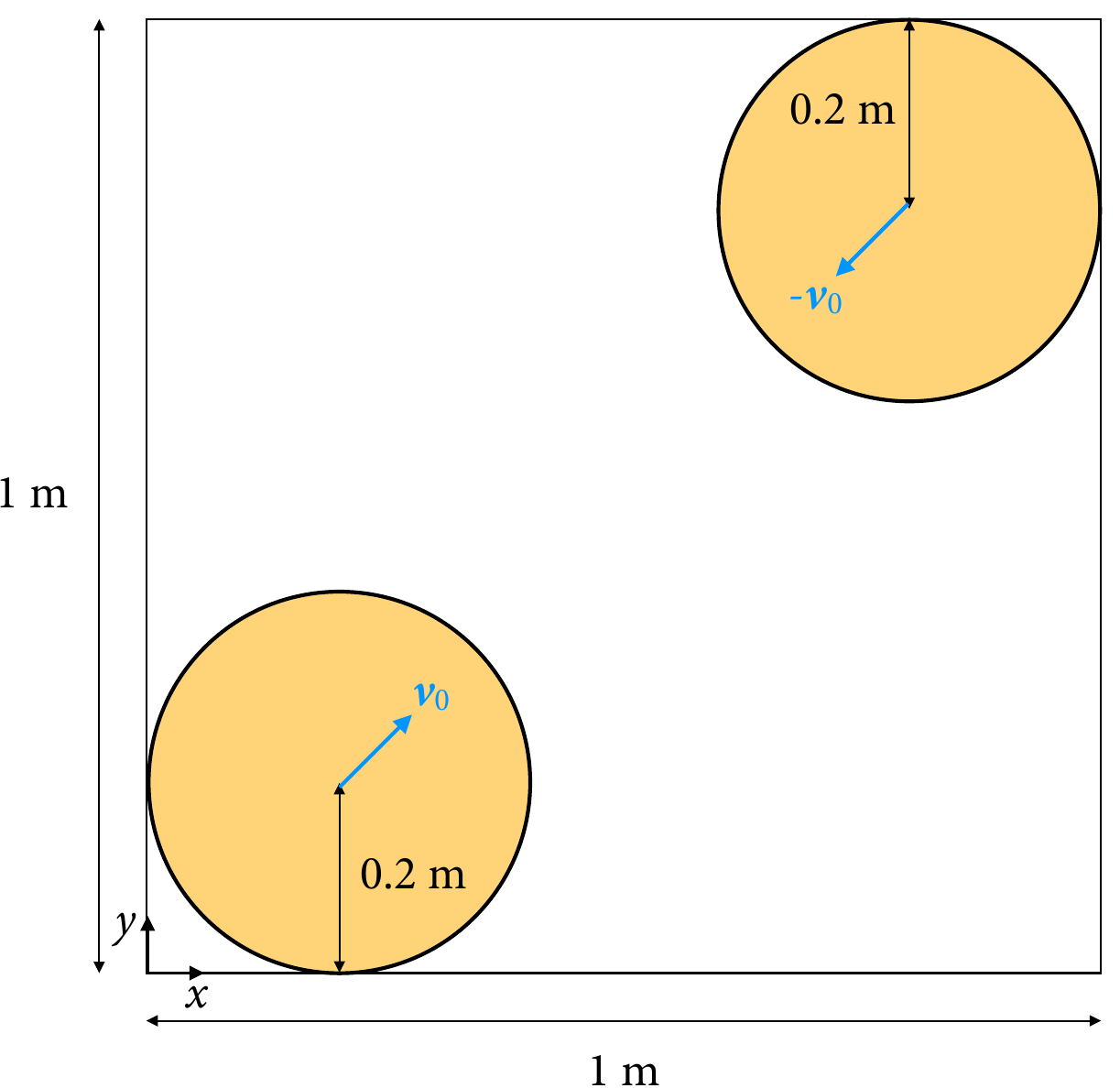}
  \caption{Setup of the impact of two poroelastic bodies problem.}
  \label{fig:impact-setup}
\end{figure}

The fourth and final example is a poromechanical extension of the impact of two elastic bodies problem which has been simulated in many papers in the MPM literature (\eg~\cite{Sulsky1994,Sulsky2004,Wang2016,Sinaie2017}).
Figure~\ref{fig:impact-setup} illustrates the setup of this problem.
Here, two circular discs of the same radius of 0.2 m are initially located at the lower left corner and the upper right corner of a 1 m square domain.
The two discs are subjected to initial velocities of the same magnitude but in opposite directions, namely $\tensor{v}_{0}=\{0.1,0.1\}$ m/s for the lower one and $-\tensor{v}_{0}=\{-0.1,-0.1\}$ m/s for the upper one.
Gravity and damping effects are neglected, and plane strain conditions are assumed for both bodies.
These conditions will make the two bodies collide at the center of the domain after a certain period of time.

Our example differs from the original problem in that here the bodies are saturated porous materials instead of solids.
Therefore, additional material parameters are required for porous materials, and they are set as follows: $\phi_{0}=0.5$, $k_{0}=10^{-14}$ m$^2$, and $\mu_{f}=10^{-6}$ kPa$\cdot$s.
The permeability is again related to porosity evolution.
For other parameters, we set the mass densities of the solid and fluid phases as $\rho_{s}=\rho_{f}=1000$ kg/m$^3$,
and the solid elasticity parameters as $E = 1000$ Pa and $\nu=0.3$, such that they are consistent with the non-dimensional parameters in the original problem.

For MPM simulation, we first introduce a background grid where the 1 m square domain in Fig.~\ref{fig:impact-setup} is discretized by mono-sized ($h=0.01$ m) square elements.
We then augment ghost elements along the boundary to ensure solvability in the very first step.
Subsequently, we generate 20,108 material points per body.
In each body, we assume that 18 material points at the top center are under zero pressure boundary conditions to allow excess pore pressure to be drained after collision.
All particles in a body are assigned the same initial velocity.
\revised{Following the original problem, we use the MPM's default contact condition, namely a no-slip condition between particles.}
\revised{Using the same time increment as the original problem for solid bodies~\cite{Sulsky1994}, namely $\Del{t}=0.001$ second}, we simulate the problem until $t=3$ seconds at which the two bodies are fully separated from the collision.
When the non-stabilized standard MPM was used, the numerical solution exhibited checkerboard oscillations from the very first step.
However, because our interest here is the stabilization of a impact problem, we suppress the stabilization term after collision for obtaining the non-stabilized standard MPM result below.

\begin{figure}[htbp]
  \centering
  \includegraphics[width=0.8\textwidth]{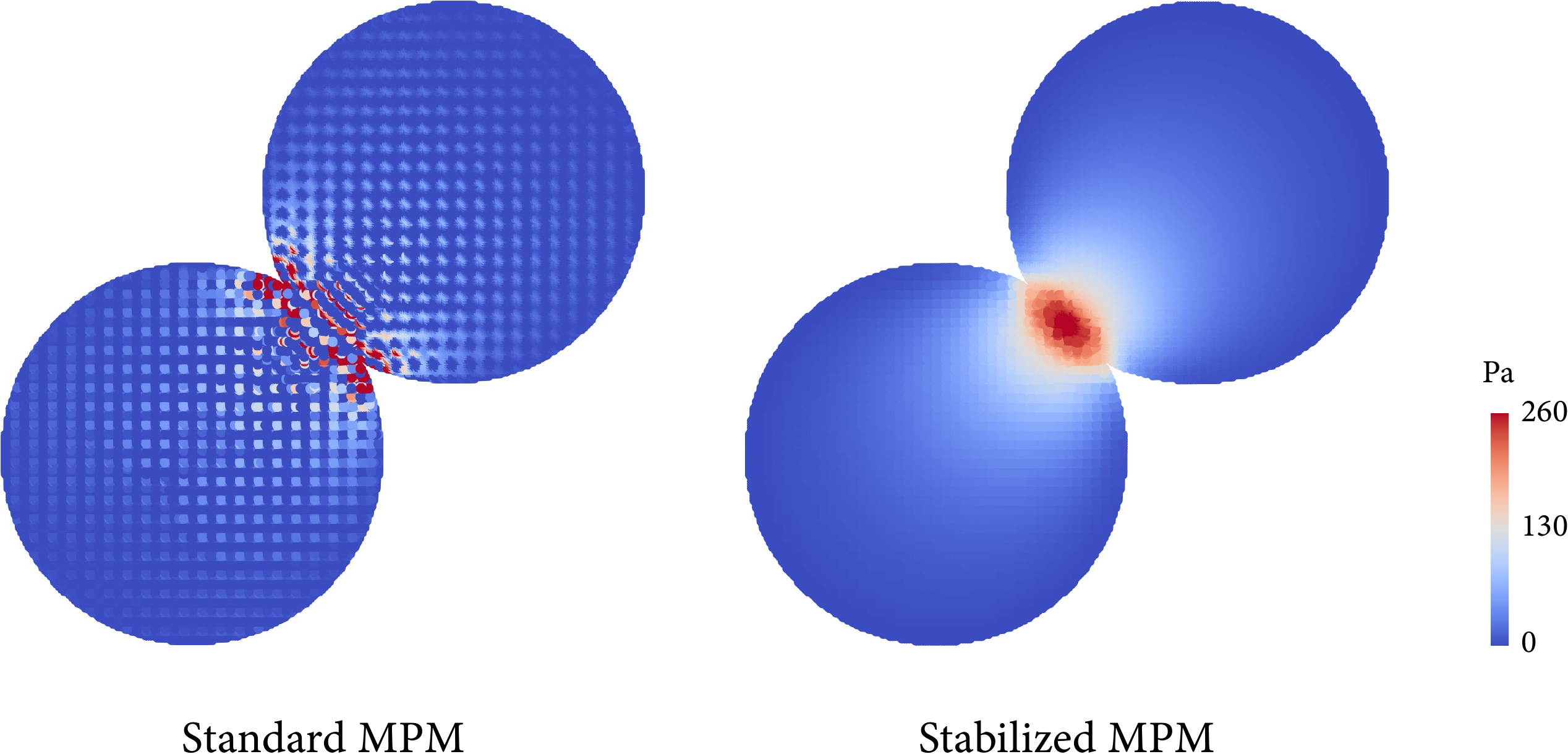}
  \caption{Excess pore pressure fields at $t=1.6$ second.}
  \label{fig:impact-collision}
\end{figure}

\begin{figure}[htbp]
  \centering
  \includegraphics[width=1.0\textwidth]{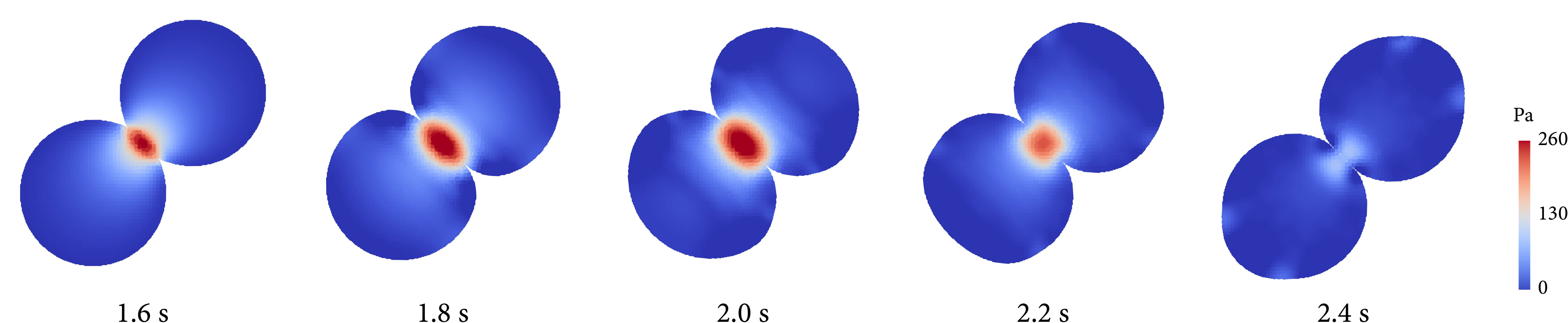}
  \caption{Excess pore pressure fields from $t=1.6$ second. All results were obtained by the stabilized MPM.}
  \label{fig:impact-contours}
\end{figure}

Figure~\ref{fig:impact-collision} compares the standard and stabilized MPM solutions of excess pore pressure fields when the two bodies collide at $t=1.6$ seconds.
Being consistent with previous results, we see that the standard MPM shows spurious pressure oscillations due to its lack of inf--sup stability.
The stabilized MPM formulation again addresses this instability remarkably well.
As far as we aware, this is the first time that a stabilized poromechanical formulation has been applied to a dynamic impact problem in the finite deformation range.
The result of this example clearly demonstrates that the stabilization method can also be useful for this type of challenging problem.
For completeness, in Fig.~\ref{fig:impact-contours} we present excess pore pressure fields obtained by the stabilized MPM from $t=1.6$ seconds.
It can be seen that the stabilized MPM provides stable pore pressure solutions throughout the collision and bouncing processes.
\revised{Note that the two bodies are subjected to tension when they are being detached from the sticky contact ($t=2.4$ seconds in Fig.~\ref{fig:impact-contours}), because here the MPM's default no-slip contact condition is used following the original problem.
It is expected that the stabilization method will also work well for other contact conditions, especially if one uses a constraint-free contact formulation such as that proposed by Fei and Choo~\cite{Fei2019}.}

% SECTION 6
% ------------------------------------------------------------------------------
\section{Closure}
\label{sec:closure}
We have introduced stabilized MPM formulations for dynamic and quasi-static poromechanics that allow for the use of standard low-order interpolation functions throughout the entire range of drainage conditions.
In these formulations, stabilization terms are augmented to the balance of mass to make equal-order mixed discretization stable even in the undrained limit.
The stabilization terms for dynamic and quasi-static MPM formulations have been commonly derived using the PPP method, but their specific expressions are different because of time-integration algorithms.
Regardless, both the dynamic and quasi-static stabilization terms can be easily implemented without change in the preexisting discretization methods.
This work has particularly used fully-implicit methods, in conjunction with the GIMP method, for both dynamic and quasi-static problems.
We have also explained how to evaluate element-wise averages in the stabilization terms when the GIMP method is used.
Through several numerical examples, we have shown that the proposed MPM formulations enable the use of standard low-order interpolation functions for undrained poromechanical problems under dynamic and quasi-static conditions.
\revised{It has also been demonstrated that the existing expressions for the stabilization parameter provide virtually identical results, suggesting that any of these expressions may be used at one's discretion.}

The key features of the stabilized MPM formulations are their efficiency and versatility.
As explained previously, the stabilization terms can be assembled using the original MPM and the GIMP method, and they are compatible with standard time-integration algorithms.
This feature allows one to simulate undrained poromechanical problems with a fairly simple modification of a standard MPM code using low-order interpolation functions.
Furthermore, previous works on mixed finite elements for poromechanics have shown that the PPP method is versatile enough to accommodate nonlinear material behavior~\cite{Sun2013a,Sun2015}, anisotropic permeability~\cite{Choo2019,Zhang2019}, and/or more field variables~\cite{Choo2015,Sun2015}.
This feature of the PPP method is expected to remain valid for the MPM as long as the relevant stabilization term is implemented in the way we described in this work.
Therefore, we believe that the proposed stabilized MPM formulations can significantly improve the efficiency of MPMs for a wide range of poromechanical problems.

% ACKNOWLEGEMENTS
\section*{Acknowledgments}
The authors are grateful to anonymous reviewers for their careful expert reviews.
This work was supported by the Research Grants Council of Hong Kong (Project 27205918) and by The University of Hong Kong (Project 201811159028).
The first author also acknowledges financial support from the Hong Kong PhD Fellowship.

% REFERENCES
\section*{References}
\bibliography{references}
% \printbibliography

\end{document}